\documentclass{amsart}

\usepackage{amsmath,amssymb,amsthm,mathtools}
\usepackage[initials]{amsrefs}
\usepackage{enumitem,microtype}
\usepackage[dvipsnames]{xcolor}
\usepackage[colorlinks=true,linkcolor=blue!55!black,citecolor=green!45!black,urlcolor=red!60!black]{hyperref}

\numberwithin{equation}{section}
\usepackage{fourier}
\usepackage[margin=1.5in]{geometry}

\newtheorem{theorem}{Theorem}[section]
\newtheorem{lemma}[theorem]{Lemma}
\newtheorem{proposition}[theorem]{Proposition}
\newtheorem{corollary}[theorem]{Corollary}
\theoremstyle{definition}
\newtheorem{definition}[theorem]{Definition}
\theoremstyle{remark}

\newcommand{\mc}[1]{\mathcal{#1}}
\newcommand{\mf}[1]{\mathsf{#1}}
\newcommand{\RR}{\mathbb R}
\newcommand{\ud}{\,\mathrm d}
\newcommand{\norm}[1]{\lVert #1\rVert}
\newcommand{\ip}[2]{\langle#1,#2\rangle}
\DeclareMathOperator{\Ent}{Ent}
\DeclareMathOperator{\TV}{TV}
\newcommand{\Ren}{\operatorname R}
\newcommand{\CLSI}{C_{\mathrm{LSI}}}
\newcommand{\CP}{C_{\mathrm P}}
\allowdisplaybreaks[1]

\usepackage{todonotes}

\title[R\'enyi mixing bounds for kinetic samplers]
      {R\'enyi mixing bounds for kinetic samplers\\
       via $L^2$ hypocoercivity}
\thanks{This work was partially supported by the National Science Foundation
under grant~\mbox{DMS-2309378}. We thank Sinho Chewi for helpful discussions.}
\author{Jianfeng Lu}
\address{Department of Mathematics, Department of Physics, and Department of Chemistry, Duke University}
\email{jianfeng@math.duke.edu}
\author{Yinchen Luo}
\address{Department of Mathematics, Duke University}
\email{yinchen.luo@duke.edu}
\date{}
\hypersetup{
  pdftitle={R\'enyi mixing bounds for kinetic samplers via L\texttwosuperior{} hypocoercivity},
  pdfauthor={Jianfeng Lu and Yinchen Luo}
}

\begin{document}

\begin{abstract}
We study convergence in R\'enyi divergence of order two for underdamped
Langevin dynamics, the bouncy particle sampler, and the Zigzag sampler,
for a target $e^{-U}$ on~$\mathbb R^d$ that satisfies a logarithmic
Sobolev inequality, with $\nabla^2U$ bounded below by a possibly negative
constant and, for the Zigzag sampler, also bounded above. For initial
laws whose R\'enyi divergence from equilibrium is at most~$R$, we show
that it decays exponentially, at a rate that depends on~$R$. With the
friction or refreshment rate chosen according to~$R$, the time needed to
reach a fixed accuracy in total variation or in R\'enyi divergence is of
order~$\sqrt{R\log R}$ at fixed dimension and target parameters,
improving the linear dependence obtained by directly converting
exponential~$L^2$ decay. The key idea is to use the logarithmic Sobolev inequality for the position target in place of the usual spatial spectral-gap estimate in the $L^2$ hypocoercivity argument. For the bouncy particle and Zigzag samplers, this
order is sharp for the R\'enyi divergence on a one-dimensional Gaussian
target, uniformly in the refreshment rate. For the
Zigzag sampler in dimension at least two, Gaussian examples show that
the mixing-time bound is sharp up to logarithmic factors in the worst
case over strongly convex targets, in both the initial divergence and the upper Hessian bound; in
particular, strong convexity alone does not give a uniform mixing bound.
\end{abstract}

\maketitle

\section{Introduction and main results}\label{sec:introduction}

We develop a general approach to establish quantitative convergence in R\'enyi divergence of order two, and apply it to three kinetic samplers: underdamped Langevin dynamics
(ULD)~\cite{ChengChatterjiBartlettJordan2018}, the bouncy particle sampler
(BPS)~\cite{BouchardCoteVollmerDoucet}, and the Zigzag sampler
(ZZ)~\cite{BierkensFearnheadRoberts}. We assume that the target measure
satisfies a logarithmic Sobolev inequality (LSI) and that its potential
has a finite Hessian lower bound, which may be negative. ZZ also requires
an upper Hessian bound, which BPS uses when it is available. For an
initial law whose R\'enyi-2 divergence from equilibrium is at most~$R$,
we prove exponential decay of this divergence at a rate depending on~$R$.
With the friction or refreshment rate chosen according to~$R$ and the
target accuracy, the time needed to reach a fixed accuracy, in total
variation or in R\'enyi-2 divergence, is of order~$\sqrt{R\log R}$ at
fixed dimension and target parameters.

Hypocoercivity estimates give exponential convergence in~$L^2$ for all
three
samplers~\cites{CaoLuWangUnderdamped,LuWangPDMP,AndrieuDurmusNueskenRoussel2021,FanLiLu26}.
Converted directly, such estimates give mixing times linear in the
initial divergence, because the $\chi^2$ divergence is exponentially
large in the R\'enyi-2 divergence. Under the LSI, we reduce the linear dependence
to~$\sqrt{R\log R}$. For ZZ in dimension~$d\ge2$, we also show that the
mixing-time bound cannot be improved by more than logarithmic factors in
the worst case over strongly convex targets, in both the initial
divergence and the upper Hessian bound, uniformly in the refreshment rate.

Throughout, we assume $U\in C^2(\RR^d)$, with
$U(x)\to\infty$ as $|x|\to\infty$ and $\int e^{-U}\ud x<\infty$. The
Hessian lower and upper bounds are denoted by $\alpha$ and $\beta$:
\begin{equation}\label{eq:hessian-bounds}
 \alpha I_d\preceq\nabla^2U(x)\preceq\beta I_d,\qquad x\in\RR^d.
\end{equation}
The lower bound $\alpha$ may be negative. The convention
$\beta=+\infty$ means that no finite upper Hessian bound is imposed. We
write $z_+=\max\{z,0\}$, so that $(-\alpha)_+$ measures the allowed
negative curvature.

The position target and the joint equilibrium measure are
\[
 \mu_x(\ud x)\propto e^{-U(x)}\ud x,\qquad
 \mu=\mu_x\otimes\mu_v,\qquad \mu_v=N(0,I_d).
\]
We assume that $\mu_x$ satisfies the LSI with constant $\CLSI>0$:
\begin{equation}\label{eq:lsi}
 \Ent_{\mu_x}(h^2)\le\frac2{\CLSI}\norm{\nabla_xh}_{L^2(\mu_x)}^2.
\end{equation}
Here $\Ent_{\mu_x}(h)=\int h\log(h/\int h\ud\mu_x)\ud\mu_x$ for
nonnegative $h$. The LSI implies the Poincar\'e inequality with the same
constant $\CLSI$. If $\alpha>0$, strong convexity gives \eqref{eq:lsi}
with $\CLSI=\alpha$ by the Bakry--\'Emery
criterion~\cite{BakryEmery1985}; see also~\cite{BGL14}*{Proposition~5.7.1}.
We keep $\CLSI$ separate from $\alpha$, since the LSI can hold when the
available Hessian lower bound is zero or negative.

ULD couples Hamiltonian transport to friction and noise in the velocity.
BPS and ZZ are piecewise deterministic Markov processes (PDMPs) whose
position moves at constant velocity between velocity jumps: bounces for
BPS, coordinate flips for ZZ, and refreshments for both. Let
$\Pi_vh(x)=\int h(x,v)\ud\mu_v(v)$ denote velocity averaging, and write
$\gamma>0$ for the friction in ULD and the refreshment rate in BPS and
ZZ. For a test function $h\in C_c^\infty(\RR^{2d})$, the three
generators are
\begin{align}
 \mc Lh
 &=(v\cdot\nabla_x-\nabla U\cdot\nabla_v)h
       +\gamma(\Delta_vh-v\cdot\nabla_vh)
 &&\text{(ULD)},\label{eq:uld}\\
 \mc Lh
 &=v\cdot\nabla_xh+(v\cdot\nabla U)_+
       [h(x,R_xv)-h(x,v)]+\gamma(\Pi_v h-h)
 &&\text{(BPS)},\label{eq:bps}\\
 \mc Lh
 &=v\cdot\nabla_xh+\sum_{i=1}^d(v_i\partial_{x_i}U)_+
       [h(x,F_iv)-h(x,v)]+\gamma(\Pi_v h-h)
 &&\text{(ZZ)}.\label{eq:zz}
\end{align}
Here $R_xv=v-2(v\cdot n)n$ reflects $v$ across the hyperplane orthogonal
to $n=\nabla U/|\nabla U|$. At critical points of $U$, choose any unit
vector for $n$; the bounce rate is zero there. The map $F_iv=v-2v_ie_i$
flips the $i$th velocity coordinate. In both PDMPs, the last term redraws
the entire velocity from $\mu_v$ at rate $\gamma$.

Under the assumptions above, the three dynamics are well posed and
preserve $\mu$. Let $\mc P_t=e^{t\mc L}$ denote the Markov semigroup of
each dynamics, and let $*$ denote the adjoint in $L^2(\mu)$. Both
$\mc P_t$ and $\mc P_t^*$ are strongly continuous contractions on
$L^2(\mu)$ that preserve positivity, constants, and integrals against
$\mu$. We write $\varrho_t$ for the law of the process at time~$t$.

For a probability measure~$\varrho\ll\mu$, the R\'enyi divergence of
order two is
\begin{equation}\label{eq:renyi-definition}
 \Ren_2(\varrho\mid\mu)=\log\int\left(\frac{\ud\varrho}{\ud\mu}\right)^2\ud\mu=\log\bigl(1+\chi^2(\varrho\mid\mu)\bigr),
\end{equation}
where $\chi^2(\varrho\mid\mu)=\norm{\ud\varrho/\ud\mu-1}_{L^2(\mu)}^2$ is
the $\chi^2$ divergence. We measure sampling accuracy by the total variation distance
\[
 \TV(\varrho,\mu)
 =\frac12\int\left|\frac{\ud\varrho}{\ud\mu}-1\right|\ud\mu.
\]
By Jensen's inequality, $\Ren_2$ bounds the Kullback--Leibler
divergence. Pinsker's inequality therefore gives
\begin{equation}\label{eq:tv-renyi}
 \TV(\varrho,\mu)^2
 \le\frac12\Ren_2(\varrho\mid\mu).
\end{equation}

\begin{theorem}[Exponential decay of R\'enyi-2 divergence]\label{thm:main}
Assume \eqref{eq:lsi} and \eqref{eq:hessian-bounds} with $\alpha\in\RR$, where
$\beta=+\infty$ is allowed for ULD and BPS and $\beta<\infty$ is required
for ZZ. Let $\varrho_t$ be the law at time $t$ of ULD, BPS, or ZZ with
friction or refreshment rate $\gamma>0$, and let $R\ge0$ satisfy
$\Ren_2(\varrho_0\mid\mu)\le R$. Then, for every $t\ge0$,
\begin{equation}\label{eq:decay}
 \Ren_2(\varrho_t\mid\mu)
 \le3e^{-\Lambda t}\Ren_2(\varrho_0\mid\mu),
\end{equation}
where
\begin{align}
 \Lambda&=\frac{2\CLSI\gamma}
 {9[\gamma^2+(-\alpha)_++3\CLSI(1+R)]}
 &&\text{for ULD};
 \label{eq:rate-uld}\\
 \Lambda&=\frac{2\CLSI\gamma}
 {900\bigl[\gamma^2+\min\{d[\CLSI(1+R)+(-\alpha)_+],\,
 \CLSI(1+R)+(-\alpha)_++\beta d\}\bigr]}
 &&\text{for BPS};
 \label{eq:rate-bps}\\
 \Lambda&=\frac{2\CLSI\gamma}
 {900[\gamma^2+\beta+(-\alpha)_++\CLSI(1+R)]}
 &&\text{for ZZ}.
 \label{eq:rate-zz}
\end{align}
\end{theorem}

The LSI is used only in the nonlinear spatial resolvent estimate,
Proposition~\ref{prop:nonlinear-resolvent}, and negative curvature enters
only through~$(-\alpha)_+$, which vanishes when~$\alpha\ge0$. The minimum
in~\eqref{eq:rate-bps} selects the better of two bounds for the BPS
cross term. One uses only the lower Hessian bound and produces the
$\beta$-free entry~$d[\CLSI(1+R)+(-\alpha)_+]$, with a dimension factor. The other
uses a finite upper Hessian bound and produces the
entry~$\CLSI(1+R)+(-\alpha)_++\beta d$. When~$\beta=+\infty$,
the minimum reduces to~$d[\CLSI(1+R)+(-\alpha)_+]$. The ZZ rate has
no dimension factor.

For $\delta>0$, the total variation mixing time
is
\[
 T_\delta^{\mathrm{TV}}=\inf\{t\ge0:\TV(\varrho_t,\mu)\le\delta\}.
\]

\begin{corollary}[R\'enyi-2 and total variation mixing]\label{cor:mixing}
Under the assumptions of Theorem~\ref{thm:main}, let $R\ge1$
satisfy $\Ren_2(\varrho_0\mid\mu)\le R$, and let
$0<\eta\le1/2$. Write $\mathfrak{L}=\log(3R/\eta)$ and choose
\begin{equation}\label{eq:time-optimized-parameters}
 \gamma^2=
 \begin{cases}
  3\CLSI(1+R/\mathfrak{L})+(-\alpha)_+,&\text{ULD},\\
  \min\{d[\CLSI(1+R/\mathfrak{L})+(-\alpha)_+],\,
        \CLSI(1+R/\mathfrak{L})+(-\alpha)_++\beta d\},&\text{BPS},\\
  \CLSI(1+R/\mathfrak{L})+\beta+(-\alpha)_+,&\text{ZZ}.
 \end{cases}
\end{equation}
Then $\Ren_2(\varrho_t\mid\mu)\le\eta$ for every time~$t$ with
\begin{equation}\label{eq:optimized-mixing-times}
 t\ge\frac{\mathfrak{L}}{\sqrt{\CLSI}}
 \begin{cases}
  60\sqrt{3(1+R/\mathfrak{L})+(-\alpha)_+/\CLSI},&\text{ULD},\\[6pt]
  6000\sqrt{\min\{d[1+R/\mathfrak{L}+(-\alpha)_+/\CLSI],\,
        1+R/\mathfrak{L}+[(-\alpha)_++\beta d]/\CLSI\}},&\text{BPS},\\[6pt]
  6000\sqrt{1+R/\mathfrak{L}+[\beta+(-\alpha)_+]/\CLSI},&\text{ZZ}.
 \end{cases}
\end{equation}
For $0<\delta\le1/2$, take $\eta=2\delta^2$, so that
$\mathfrak{L}=\log[3R/(2\delta^2)]$. Then the total variation mixing time
$T_\delta^{\mathrm{TV}}$ is at most the right-hand side
of~\eqref{eq:optimized-mixing-times}.
\end{corollary}

The total variation bound follows from the R\'enyi-2 bound
by~\eqref{eq:tv-renyi}. By data processing,
both bounds also hold for the position marginal. The numerical constants
are not optimized.

The initial divergence can be large in high dimension. Suppose that
$0<\alpha\le\beta<\infty$, and let~$x_\star$ minimize~$U$. For the
Gaussian cold start $\varrho_0=N(x_\star,\beta^{-1}I_d)\otimes\mu_v$, the
position density relative to~$\mu_x$ is at
most~$(\beta/\alpha)^{d/2}$~\cite{DwivediChenWainwrightYu2019}, so
$\Ren_2(\varrho_0\mid\mu)\le\frac d2\log(\beta/\alpha)$.
Table~\ref{tab:intro-leading} compares the bounds of
Corollary~\ref{cor:mixing} for $R\asymp d$ up to logarithmic
factors, against the direct conversion of exponential~$L^2$ decay and
against bounds obtained by other methods. For this comparison, all rows assume that the target satisfies the LSI and that~$U$ is convex. The table applies to initial laws with $\Ren_2(\varrho_0\mid\mu)\le R$; the Gaussian cold start above motivates the choice of scale. Our results also allow $\alpha<0$, with the full dependence on~$\alpha$ given in Corollary~\ref{cor:mixing}. Section~\ref{sec:related-work}
discusses these methods and compares the resulting gradient-query
complexities with earlier bounds.

\begin{table}[htbp]
\centering
\renewcommand{\arraystretch}{1.3}
\begin{tabular}{@{}lll@{}}
\hline
Sampler & Additional assumptions & Mixing time $T_\delta^{\mathrm{TV}}$ for $R\asymp d$\\
\hline
\multicolumn{3}{@{}l}{\emph{Direct conversion of $L^2$ hypocoercivity}}\\
ULD~\cites{CaoLuWangUnderdamped,FanLiLu26} & none & $\CLSI^{-1/2}d$\\
BPS~\cite{LuWangPDMP} & none & $\CLSI^{-1/2}d^{3/2}$\\
ZZ~\cite{LuWangPDMP} & $\nabla^2U\preceq\beta I_d$ & $\CLSI^{-1}\sqrt\beta\,d$\\
\hline
\multicolumn{3}{@{}l}{\emph{This work, Corollary~\ref{cor:mixing}}}\\
ULD & none & $\CLSI^{-1/2}\sqrt d$\\
BPS & none & $\CLSI^{-1/2}d$\\
BPS & $\nabla^2U\preceq\beta I_d$ & $\CLSI^{-1}\sqrt{\beta d}$\\
ZZ & $\nabla^2U\preceq\beta I_d$ & $\CLSI^{-1/2}\sqrt d+\CLSI^{-1}\sqrt\beta$\\
\hline
\multicolumn{3}{@{}l}{\emph{Entropy hypocoercivity}}\\
ULD~\cite{LiLuHypercontractivity} & none & $\CLSI^{-1/2}$\\
BPS~\cite{LuLuoEntropy} & $0<\alpha\le\nabla^2U\le\beta$ ($\CLSI=\alpha$); bounded initial density & $\CLSI^{-1/2}\sqrt d$\\
\hline
\multicolumn{3}{@{}l}{\emph{Lower bound for ZZ, Corollary~\ref{cor:zz-sharpness}}}\\
ZZ & $\nabla^2U\preceq\beta I_d$ & $\CLSI^{-1/2}\sqrt d+\CLSI^{-1}\sqrt\beta$\\
\hline
\end{tabular}
\bigskip 
\caption{Bounds on the total variation mixing time for initial laws with $\Ren_2(\varrho_0\mid\mu)\le R$, at the scale $R\asymp d$, up to universal constants and logarithmic factors in the relevant parameters. The last block is a lower bound, and the others are upper bounds. All rows assume the LSI~\eqref{eq:lsi} and $\nabla^2U\succeq0$, and technical regularity conditions on~$U$ in the cited results are omitted.  The first block converts the exponential~$L^2$ decay of~\cites{CaoLuWangUnderdamped,FanLiLu26} for ULD and of~\cite{LuWangPDMP}*{Theorem~1} for BPS and ZZ, as explained after the table. The second block follows from Corollary~\ref{cor:mixing}, which also gives the full dependence on~$R$ and the case $\alpha<0$. Under an upper Hessian bound, BPS satisfies the smaller of its two entries, and the entry $\CLSI^{-1}\sqrt{\beta d}$ absorbs a term $\CLSI^{-1/2}\sqrt d$ because $\CLSI\le\beta$. The third block is taken from~\cite{LiLuHypercontractivity}*{Corollary~5.5} and~\cite{LuLuoEntropy}*{Theorem~1.1}; in its BPS row, the initial density relative to~$\mu$ is at most~$e^R$, which implies $\Ren_2(\varrho_0\mid\mu)\le R$. By Corollary~\ref{cor:zz-sharpness}, some Gaussian target satisfying the assumptions of the last row, with an initial law satisfying $\Ren_2(\varrho_0\mid\mu)\le R$, has at least the mixing time shown there for every refreshment rate. Hence, for $d\ge2$, the ZZ entry of the second block is sharp up to the suppressed factors.}
\label{tab:intro-leading}
\end{table}

As $R\to\infty$ with the other parameters fixed, the bounds of
Corollary~\ref{cor:mixing} are of order~$\CLSI^{-1/2}\sqrt{R\log R}$,
with an extra factor~$\sqrt d$ for BPS without an upper Hessian bound;
negative curvature does not change this order. By comparison, an $L^2$ estimate
$\chi^2(\varrho_t\mid\mu)\le Ce^{-\lambda t}\chi^2(\varrho_0\mid\mu)$
gives
\[
 \Ren_2(\varrho_t\mid\mu)\le\log\left[1+C e^{-\lambda t}\bigl(e^{\Ren_2(\varrho_0\mid\mu)}-1\bigr)\right],
\]
which guarantees a fixed accuracy only after a time of
order~$(R+\log C)/\lambda$. The rate~$\lambda$ does not depend on the
initial law, so optimizing it over the friction or refreshment rate
still leaves this time linear in~$R$. The first block of
Table~\ref{tab:intro-leading} records this time, at $R\asymp d$, for the
optimized~$L^2$ rates
of~\cites{CaoLuWangUnderdamped,FanLiLu26,LuWangPDMP}, which are of
order~$\sqrt{\CLSI}$ for ULD, $\sqrt{\CLSI/d}$ for BPS, and
$\CLSI/\sqrt\beta$ for ZZ, with the Poincar\'e constant replaced
by~$\CLSI$. These $L^2$ results use the Poincar\'e inequality in place
of the LSI. The PDMP rates in~\cite{LuWangPDMP} are stated under a
compactness condition on~$\mu_x$ that does not enter them. A single
application of
Theorem~\ref{thm:main}, with~$\gamma$ optimized for~$R$, already gives
order~$\CLSI^{-1/2}\sqrt R\log R$. The proof of
Corollary~\ref{cor:mixing} applies the theorem in stages, halving the
divergence bound at each stage, and removes a further
factor~$\sqrt{\log R}$.

Following Fan, Li, and Lu~\cite{FanLiLu26}, we add to the quadratic
energy a small multiple of a shifted corrector, whose time derivative
dissipates the spatial component of the density through a macroscopic
dissipation. A Poincar\'e inequality bounds this dissipation below by a
multiple of the squared norm, which gives exponential $\chi^2$ decay and
hence mixing times linear in~$R$. Our main new estimate is a nonlinear
lower bound on the macroscopic dissipation, obtained by applying the LSI
to a smoothed position density; it controls the position marginal on the
logarithmic scale of the R\'enyi-2 divergence. Choosing the shift to be $\CLSI(1+R)$ increases the uniform denominator bound by only a constant factor compared with zero shift. This larger shift admits a larger corrector weight, which raises the rate from
order~$\sqrt{\CLSI}/R$ to order~$\sqrt{\CLSI/R}$. The velocity dissipation
controls the gap between the joint density and its position marginal.

\begin{theorem}[Sharpness for the R\'enyi-2 divergence]\label{thm:initial-divergence-sharpness}
Let $d=1$ and $\alpha>0$, set $U(x)=\alpha x^2/2$ on~$\RR$, and write
$\mu=N(0,\alpha^{-1})\otimes N(0,1)$. For every~$R\ge e^8$, there
exists a probability law~$\varrho_0$ with~$\Ren_2(\varrho_0\mid\mu)=R$
such that, for BPS and ZZ with every refreshment rate~$\gamma>0$,
\begin{equation}\label{eq:initial-divergence-lower-bound}
 \Ren_2(\varrho_t\mid\mu)>\frac12
 \qquad\text{for all }0\le t<\frac1{16\sqrt\alpha}\sqrt{R\log R}.
\end{equation}
\end{theorem}

With~$\CLSI=\alpha$ and~$d=1$, Corollary~\ref{cor:mixing} with
$\eta=1/2$ gives $\Ren_2(\varrho_t\mid\mu)\le1/2$ after a
time of order~$\alpha^{-1/2}\sqrt{R\log R}$, so this order is sharp.
Theorem~\ref{thm:initial-divergence-sharpness} also limits the decay
rate. If $\Ren_2(\varrho_t\mid\mu)\le3e^{-\Lambda t}R$ for all $t\ge0$, as
in~\eqref{eq:decay}, then evaluating at $t=\log(6R)/\Lambda$ gives
\[
 \Lambda\le\frac{16\sqrt\alpha\log(6R)}{\sqrt{R\log R}}.
\]
For this target, optimizing~$\gamma$ in~\eqref{eq:rate-bps}
and~\eqref{eq:rate-zz} gives rates of order~$\sqrt{\alpha/R}$. Hence no
refreshment rate improves the optimized rate of Theorem~\ref{thm:main} by
more than a factor of order~$\sqrt{\log R}$. The obstruction is specific
to velocity jumps. Bounces and flips preserve the speed~$|v|$, so slow
refreshment leaves the initial speed distribution in place, while
frequent refreshment slows the return from a distant initial position.
The proof is in Appendix~\ref{app:initial-divergence-sharpness}.

For ZZ in dimension~$d\ge2$,\footnote{For dimension $d = 1$, ZZ coincides with BPS.} the bound of Corollary~\ref{cor:mixing} is
sharp up to logarithmic factors in the worst case over strongly convex
targets and initial laws with $\Ren_2(\varrho_0\mid\mu)\le R$, in both the
initial divergence and the upper Hessian bound. This concerns the whole
class of such initial laws, not the Gaussian cold start specifically.
Appendix~\ref{sec:zz-obstruction} proves the sharpness with two Gaussian examples,
each uniform in the refreshment rate. 

In the first,
Theorem~\ref{thm:zz-finite-speed}, the target is isotropic and the
initial law is displaced by~$\sqrt{R/\CLSI}$ along a coordinate axis.
Flips preserve the absolute value of each velocity coordinate, so each
position coordinate moves at bounded mean speed, and the total variation
distance stays above~$1/4$ for a
time of order~$\sqrt{R/\CLSI}$.

The second result, Theorem~\ref{thm:zz-obstruction}, addresses the
dependence on~$\beta$. It exploits the coordinatewise design of ZZ, in
which each event reverses a single velocity coordinate at a rate set by
the corresponding partial derivative of~$U$. The targets are
two-dimensional Gaussians with~$\CLSI=1$, increasingly confined in the
direction of~$x_1-x_2$, which is oblique to the coordinate axes. Both
partial derivatives of~$U$ then grow with the confinement~$\beta$. The
resulting frequent flips randomize the velocity component along the
valley~$\{x_1=x_2\}$ before the particle can travel along it. The
initial laws have the same density relative to the target for
every~$\beta$, so their initial divergences do not depend on~$\beta$.
Nevertheless, the total variation distance and the R\'enyi-2 divergence
stay above half their initial values for a time of
order~$\sqrt\beta/\log(e\beta)$, uniformly in the refreshment rate. This
matches the term~$\CLSI^{-1}\sqrt\beta$
in~\eqref{eq:optimized-mixing-times} up to the logarithm.
Corollary~\ref{cor:zz-sharpness} combines the two examples into a lower
bound on the ZZ mixing time, which shows that the ZZ entry of
Corollary~\ref{cor:mixing} in Table~\ref{tab:intro-leading} cannot be
improved by more than logarithmic factors. Since this lower bound is
unbounded in~$\beta$, strong convexity alone gives no uniform ZZ bound.
BPS avoids this obstruction: for this family, the $\beta$-free entry
of~\eqref{eq:rate-bps} gives a BPS bound uniform in the confinement.
Heuristically, a bounce reflects the velocity across~$\nabla U$ instead
of reversing a single coordinate. When the confinement is strong,
$\nabla U$ is nearly orthogonal to the valley on most of the target, so
a bounce nearly preserves the velocity component along the valley.

\subsection{Related works}\label{sec:related-work}

ULD underlies several discretized sampling algorithms; see, for
example,~\cites{ChengChatterjiBartlettJordan2018,AltschulerChewiZhang2026}
and the book~\cite{ChewiLogConcaveSampling}. The event-driven bounce
mechanism of BPS originated with Peters and de With~\cite{PetersDeWith2012}
and was developed into the BPS framework by Bouchard-C\^ot\'e, Vollmer, and
Doucet~\cite{BouchardCoteVollmerDoucet}. A one-dimensional Zigzag
process arose as a scaling limit of lifted Metropolis--Hastings chains for the
Curie--Weiss model~\cite{BierkensRoberts2017}, and Bierkens, Fearnhead,
and Roberts~\cite{BierkensFearnheadRoberts} developed the Zigzag process as
a general-purpose sampler. Both belong to the broader class of PDMP Monte
Carlo methods surveyed in~\cite{FearnheadEtAl2018}. Harris-type arguments
establish geometric ergodicity of BPS and ZZ under tail and
curvature conditions~\cites{DeligiannidisBouchardCoteDoucet2019,DurmusGuillinMonmarche2020,BierkensRobertsZitt2019}.
Bierkens, Kamatani, and Roberts~\cite{BierkensKamataniRoberts2025}
study the contrasting behavior of BPS and ZZ on anisotropic Gaussian
targets through scaling limits.

In kinetic samplers, the position follows a velocity that is randomized
only gradually, so on short time scales it moves ballistically rather
than diffusively. This persistent motion can accelerate convergence. For
convex potentials, ULD with optimized friction converges in~$L^2$ at a
rate of order~$\sqrt{\CP}$, where $\CP$ is the Poincar\'e constant
of~$\mu_x$, whereas the overdamped Langevin diffusion converges at
rate~$\CP$~\cites{CaoLuWangUnderdamped,FanLiLu26}. Hypocoercivity
explains how Hamiltonian transport turns velocity dissipation into
convergence in position. Villani~\cite{Villani2009} developed a general
theory of hypocoercivity for kinetic equations. Dolbeault, Mouhot, and
Schmeiser~\cites{DMS09,DMS15} introduced the modified~$L^2$ method.
Grothaus and Stilgenbauer~\cite{GS16} gave a rigorous operator
treatment for Langevin dynamics, and Andrieu, Durmus, N\"usken, and
Roussel~\cite{AndrieuDurmusNueskenRoussel2021} developed modified~$L^2$
hypocoercivity for a broad class of velocity-jump samplers. Fan,
Li, and Lu~\cite{FanLiLu26}*{Theorem~1 and Corollary~2} developed the
shifted corrector under a Poincar\'e inequality of~$\mu_x$ and a Hessian
lower bound. At zero shift, optimizing the friction recovers the
sharp square-root dependence of the rate on the Poincar\'e constant for
convex potentials. Their framework is the starting point for our
argument.

A complementary route to~$L^2$ hypocoercivity uses space-time
Poincar\'e inequalities. Cao, Lu,
and Wang~\cite{CaoLuWangUnderdamped} obtained explicit~$L^2$ rates for ULD.
Lu and Wang~\cite{LuWangPDMP}*{Theorem~1} obtained such rates for
randomized Hamiltonian Monte Carlo, BPS, and ZZ. From the lifting
viewpoint, Eberle and L\"orler~\cite{EberleLorlerLifts} interpret BPS
and ZZ as nonreversible lifts and show that any second-order lift can
improve the relaxation time by at most a square-root factor. All these
results give quantitative $L^2$ convergence; the improvement in our
mixing bounds is relative to the direct conversion of such
exponential~$L^2$ estimates, under the additional assumption that~$\mu_x$
satisfies an LSI.

For ULD, Lu~\cite{LuUnderdampedEntropy} established sharp hypocoercive
entropy decay using the Wasserstein entropy current as a corrector. Li and
Lu~\cite{LiLuHypercontractivity}*{Corollary~5.5} developed a space-time LSI
and hypocoercive hypercontractivity. For convex potentials, under an LSI of~$\mu_x$ and a
Hessian-growth condition, they obtain exponential decay
of every R\'enyi order~$q>1$, with constants independent of the initial
law. At fixed accuracy, their mixing-time bound grows logarithmically with
the initial divergence, better than our~$\sqrt{R\log R}$ dependence
where both results apply. Our argument also applies to BPS
and~ZZ, for which logarithmic dependence on the initial divergence is
impossible, uniformly in the refreshment rate.
Theorem~\ref{thm:initial-divergence-sharpness} rules it out for the time
at which the R\'enyi-2 divergence falls below~$1/2$.
Theorem~\ref{thm:zz-finite-speed} rules it out for the total variation
mixing time of ZZ in every dimension and of BPS in dimension one.

For PDMPs, Monmarch\'e and Wang~\cite{MonmarcheWangEntropy}*{Theorem~2}
show that BPS and ZZ admit no uniform multiplicative contraction of
the KL divergence over the initial laws whose KL divergence lies below a
fixed positive level, even for a standard Gaussian target. Their theorem
covers standard Gaussian velocities. Our estimates require a finite
bound~$R$ on the initial R\'enyi-2 divergence, equivalently an~$L^2$
relative density, and allow the rate to depend on~$R$. Since a KL bound
does not control~$\Ren_2$, this stronger condition on the initial law
is consistent with their obstruction.

In a companion paper~\cite{LuLuoEntropy}, we 
study relative entropy decay for BPS with strongly convex and smooth
potentials, using the Wasserstein entropy current of
Lu~\cite{LuUnderdampedEntropy} as a corrector. For initial laws with
density bounded by~$M$ relative to~$\mu$, we obtain an exponentially
decaying entropy estimate with an additive remainder, at a rate depending
logarithmically on~$M$ and on the size of the remainder. The density
bound assumption gives $\Ren_2(\varrho_0\mid\mu)\le\log M$, so this condition on the
initial law is also consistent with the obstruction of Monmarch\'e and
Wang.

Where both results apply, the companion bound has a better dependence on
the dimension and the condition number; see the BPS row of the entropy
hypocoercivity block in Table~\ref{tab:intro-leading}. The dimension enters its rate additively,
alongside~$\log M$, and the condition number only logarithmically.
In~\eqref{eq:rate-bps}, by contrast, the dimension multiplies the
divergence term or enters through~$\beta d$. For the Gaussian cold start,
this difference appears in the query complexities discussed below. In
the initial divergence alone, with $R=\log M$, the companion bound is
larger by a factor~$\sqrt{\log R}$, and it requires a bounded initial
density rather than a finite R\'enyi-2 divergence. The R\'enyi-2 route
taken here is thus less sharp for BPS in the dimension and the condition
number, but it rests on a different argument. It upgrades the $L^2$ hypocoercivity estimate through a single
nonlinear spatial resolvent estimate,
Proposition~\ref{prop:nonlinear-resolvent}. It therefore covers ULD, ZZ, and nonconvex
potentials satisfying an LSI, and it should apply to other dynamics for
which the modified $L^2$ estimate is available.

The size of the initial divergence is a recurring issue in the complexity
of sampling. Warm-start analyses of the Metropolis-adjusted Langevin
algorithm assume a bound on it~\cites{DwivediChenWainwrightYu2019,WuSchmidlerChen2022},
and Altschuler and Chewi~\cite{AltschulerChewi2024} construct R\'enyi warm
starts algorithmically using discretized ULD. For ZZ, Lu and
Wang~\cite{LuWangZigzag} bound the query complexity under a warm-start
assumption. Our bounds concern mixing in continuous time. Lu and
Luo~\cite{LuLuo} simulate BPS and ZZ exactly by windowed thinning,
with an expected number of gradient queries proportional to the time
horizon. Combining their analysis with Corollary~\ref{cor:mixing}, for
$0<\alpha\le\beta<\infty$ and~$\kappa=\beta/\alpha$, the Gaussian cold
start reaches total variation accuracy~$\delta$ with
$\widetilde O(\min\{\kappa^{1/2}d^{3/2},\kappa d\})$ gradient queries for
BPS and $\widetilde O(\kappa^{1/2}d^{3/4}+\kappa d^{1/4})$ full-gradient
equivalents for ZZ, counting~$d$ coordinate-partial queries as one.
Here~$\widetilde O$ suppresses logarithmic factors in~$d$, $\kappa$,
and~$1/\delta$. These improve the bounds $\widetilde O(\kappa^{1/2}d^2)$
and $\widetilde O(\kappa d^{5/4})$ obtained in~\cite{LuLuo}. For BPS
under the same assumptions, the entropy estimate of the companion
paper~\cite{LuLuoEntropy} gives the smaller bound
$\widetilde O(\kappa^{1/2}d)$.

\subsection*{Organization}
Section~\ref{sec:preliminaries} sets up the relative densities and recalls
the shifted corrector and the modified energy.
Section~\ref{sec:dissipation} establishes the dissipation of the modified
energy and the nonlinear spatial resolvent estimate, and ends with the
sampler-specific calculations. Section~\ref{sec:closure} proves
Theorem~\ref{thm:main} and Corollary~\ref{cor:mixing}.
Appendix~\ref{app:initial-divergence-sharpness} proves
Theorem~\ref{thm:initial-divergence-sharpness}, and
Appendix~\ref{sec:zz-obstruction} proves the lower bounds that show the
ZZ bound to be sharp up to logarithmic factors.

\subsection*{AI Usage} During our project to improve the cold-start query complexity estimate of~\cite{LuLuo}, which led to the entropy estimate in~\cite{LuLuoEntropy}, our interaction with ChatGPT 6 Astra led to the observation that a R\'enyi-2 contraction estimate could improve the query complexity estimate of BPS. After checking the AI output, whose argument was rather opaque, the authors realized the underlying mechanism: the $L^2$ hypocoercivity estimate can be upgraded by applying the LSI to the position marginal, rather than simply using coercivity. The authors then wrote the simplified proof based on their understanding, extended it to other kinetic sampling processes, and prepared the manuscript. The authors take full responsibility for its mathematical contents. Although we do not intend to submit this manuscript for publication in journals, we believe the approach of establishing better R\'enyi-2 convergence rate via $L^2$ hypocoercivity might be of general interest, and thus decide to make it available on arXiv.

\section{Preliminaries}\label{sec:preliminaries}

Unless indicated otherwise, norms and inner products are in $L^2(\mu)$;
on functions of~$x$ alone they coincide with those in $L^2(\mu_x)$. For
operators acting only on position, the adjoint is taken in
$L^2(\mu_x)$. We work with real functions and use Hilbert--Schmidt norms
for matrices. The subscript $0$ denotes the zero-mean subspace:
\begin{equation}\label{eq:zero-mean-space}
 L^2_0(\mu):=\left\{f\in L^2(\mu):\int f\ud\mu=0\right\}.
\end{equation}
We use the same convention for $L^2_0(\mu_x)$.

We identify $\Pi_v h$ with its lift to phase space, so that $\Pi_v$ is the
orthogonal projection onto functions of~$x$ alone in $L^2(\mu)$. Test
functions are elements of $C_c^\infty(\RR^{2d})$, or of
$C_c^\infty(\RR^d)$ for functions of~$x$ alone. On functions of~$x$
alone, all three generators reduce to transport:
\[
 \mc L\Pi_v=(v\cdot\nabla_x)\Pi_v.
\]
The overdamped generator $\mc L_{\mathrm o}=\Delta_x-\nabla U\cdot\nabla_x$
acts on functions of~$x$ alone. We realize $-\mc L_{\mathrm o}$ as the
nonnegative self-adjoint operator associated with the closed form
$h\mapsto\norm{\nabla_xh}^2$.

\subsection{Relative densities}\label{sec:relative-densities}

For the law $\varrho_t$ of the process at time~$t$, define
\begin{equation}\label{eq:density}
 f_t=\frac{\ud\varrho_t}{\ud\mu}-1.
\end{equation}
We decompose it into its spatial component and velocity fluctuation:
\begin{equation}\label{eq:centered-density}
 f_{S,t}=\Pi_v f_t,\qquad f_{F,t}=(I-\Pi_v)f_t.
\end{equation}
Thus $f_t=f_{S,t}+f_{F,t}$ and $\Pi_v f_{F,t}=0$.
The functions $1+f_t$ and $1+f_{S,t}$ are the relative densities of
$\varrho_t$ with respect to $\mu$ and of its position marginal with
respect to $\mu_x$. In particular,
\[
 f_t\ge-1,\qquad f_{S,t}\ge-1,\qquad
 \int f_t\ud\mu=\int f_{S,t}\ud\mu_x=0.
\]
When no confusion arises, we suppress the time subscript and write
$f,f_S,f_F$.
For every test function $h$, invariance of $\mu$ gives
$\int\mc Lh\ud\mu=0$, and hence
\begin{equation}\label{eq:weak-evolution}
 \frac{\ud}{\ud t}\ip{f_t}{h}
 =\int\mc Lh\ud\varrho_t
 =\ip{1+f_t}{\mc Lh}
 =\ip{f_t}{\mc Lh}.
\end{equation}
Consequently,
\begin{equation}\label{eq:density-evolution}
 \partial_tf_t=\mc L^*f_t.
\end{equation}
This equation is understood in the weak sense; if $f_0\in D(\mc L^*)$, it holds in $L^2(\mu)$.
The zero-mean property and orthogonality also imply
\begin{equation}\label{eq:centered-norms}
 \begin{gathered}
 \norm{1+f_t}^2=1+\norm{f_t}^2,\qquad
 \norm{1+f_{S,t}}^2=1+\norm{f_{S,t}}^2,\\
 \norm{f_t}^2=\norm{f_{S,t}}^2+\norm{f_{F,t}}^2,\qquad
 \Ren_2(\varrho_t\mid\mu)=\log(1+\norm{f_t}^2).
 \end{gathered}
\end{equation}

\subsection{The shifted corrector and modified energy}
\label{sec:modified-energy}

The quadratic energy $\frac12\norm f^2$ dissipates only the velocity
fluctuation~$f_F$ (Lemma~\ref{lem:velocity-dissipation}). A spatial
perturbation is not stationary, however: transport turns the spatial
gradient of~$f_S$ into a velocity fluctuation, which the velocity
dynamics then dissipates. The modified $L^2$ method of Dolbeault, Mouhot,
and Schmeiser~\cites{DMS09,DMS15}, which builds on ideas of
H\'erau~\cite{Herau2006}, quantifies this mechanism. It adds
to~$\frac12\norm f^2$ a small multiple of a corrector, designed so that
its time derivative dissipates~$f_S$. We use the shifted corrector of
Fan, Li, and Lu~\cite{FanLiLu26}*{Section~2}.

\begin{definition}[Shifted corrector and modified energy]
\label{def:modified-energy}
Let $a>0$ and $\epsilon>0$. The shifted corrector is
\begin{equation}\label{eq:corrector}
 \mc A_a=(a-\mc L_{\mathrm o})^{-1}
          \bigl((v\cdot\nabla_x)\Pi_v\bigr)^*.
\end{equation}
The corrector is defined first on test functions and extended by continuity to
$L^2_0(\mu)$; see Lemma~\ref{lem:corrector-bounds}.
For $f\in L^2_0(\mu)$, the modified energy is
\begin{equation}\label{eq:modified-energy}
 \mf L_{a,\epsilon}(f)
 :=\frac12\norm f^2
     +\epsilon\ip{\mc A_af}{f}.
\end{equation}
For the law $\varrho_t$, we apply this energy to
$f_t=\ud\varrho_t/\ud\mu-1$ from \eqref{eq:density}.
\end{definition}

The correction enters \eqref{eq:modified-energy} with the sign opposite
to that in~\cite{FanLiLu26}. They differentiate the energy along the
backward equation $\partial_tf=\mc Lf$, whereas we use the forward
density equation~\eqref{eq:density-evolution}, in which the transport
term of the generator changes sign.

The shift $a$ and the correction weight $\epsilon$ are auxiliary proof
parameters, distinct from the Hessian bounds $\alpha,\beta$, the LSI
constant $\CLSI$, and the friction or refreshment rate $\gamma$. Both
are constant in time. The weight is chosen in
Proposition~\ref{prop:common-dissipation}, and the shift in the proofs of
Section~\ref{sec:closure}.

We record operator bounds for the corrector. They also justify the
differentiation of the correction in the proof of
Proposition~\ref{prop:common-dissipation}.

\begin{lemma}[Corrector identities and operator bounds]
\label{lem:corrector-bounds}
For every $a>0$, the corrector in Definition~\ref{def:modified-energy}
has a bounded extension to $L^2_0(\mu)$ satisfying
\begin{equation}\label{eq:operator}
 \begin{gathered}
 \Pi_v \mc A_a=\mc A_a,\qquad
 \mc A_a\Pi_v=0,\\
 \norm{\mc A_a}\le\frac1{2\sqrt a},\qquad
 \norm{(v\cdot\nabla_x)\mc A_a}\le1.
 \end{gathered}
\end{equation}
Here the norms in \eqref{eq:operator} are operator norms on
$L^2_0(\mu)$.
Moreover, for every $f\in L^2_0(\mu)$,
\begin{equation}\label{eq:corrector-adjoint}
 \mc A_a^*f=(v\cdot\nabla_x)(a-\mc L_{\mathrm o})^{-1}\Pi_v f.
\end{equation}
\end{lemma}
The identities and bounds are proved in Fan, Li, and
Lu~\cite{FanLiLu26}*{Section~3 and Lemma~1}.

For a sufficiently small correction weight, the modified energy is
comparable to the squared $L^2$ norm.

\begin{lemma}[Energy equivalence]\label{lem:energy-comparison}
For every $a>0$, $0<\epsilon\le\sqrt a$, and $f\in L^2_0(\mu)$,
\begin{equation}\label{eq:equivalence}
 \frac14\norm f^2\le\mf L_{a,\epsilon}(f)
 \le\frac34\norm f^2.
\end{equation}
\end{lemma}

\begin{proof}
Since $\Pi_v f_F=0$ and $f_S$ is spatial,
\[
 \ip{f_F}{f_S}=\ip{\Pi_v f_F}{f_S}=0,\qquad
 \norm f^2=\norm{f_S}^2+\norm{f_F}^2.
\]
Also $\mc A_a f=\mc A_a f_F$ is spatial, so
\begin{align*}
 \bigl|\epsilon\ip{\mc A_a f}{f}\bigr|
 &=\bigl|\epsilon\ip{\mc A_a f_F}{f_S}\bigr|\\
 &\le\frac{\epsilon}{2\sqrt a}\norm{f_F}\norm{f_S}\\
 &\le\frac{\epsilon}{4\sqrt a}
       (\norm{f_F}^2+\norm{f_S}^2)
 \le\frac14\norm f^2.
\end{align*}
The last inequality uses $\epsilon\le\sqrt a$. Hence
\eqref{eq:equivalence} follows.
\end{proof}

\section{Dissipation estimates}\label{sec:dissipation}

We establish the two estimates that yield decay of the R\'enyi-2
divergence. The first bounds the time derivative of the modified energy
by a macroscopic dissipation and the velocity dissipation. The second
bounds the macroscopic dissipation below through the LSI; it is the only
place where the LSI enters.

For $f_S\in L^2_0(\mu_x)$, define the macroscopic dissipation
\begin{equation}\label{eq:Q}
 \mc D_a(f_S)
 =\ip{f_S}{(-\mc L_{\mathrm o})(a-\mc L_{\mathrm o})^{-1}f_S}.
\end{equation}
The name follows the macroscopic coercivity condition
of~\cite{DMS15}*{Assumption~(H2)} and the terminology
of~\cite{FanLiLu26}.

\begin{proposition}[Dissipation of the modified energy]\label{prop:common-dissipation}
Assume \eqref{eq:hessian-bounds}, with $\beta<\infty$ for ZZ.
Let $a>0$, and let the corrector weight $\epsilon$ be given by the formulas below. Let $\varrho_t$ be the law of the process at time~$t$, and let $f_t=\ud\varrho_t/\ud\mu-1$ as in~\eqref{eq:density}. If $f_0\in D(\mc L^*)$, then the modified energy satisfies, for all $t\ge0$,
\begin{equation}\label{eq:common-dissipation}
 \frac{\ud}{\ud t}\mf L_{a,\epsilon}(f_t)
 \le-\frac\epsilon2\mc D_a(f_S)
    -\frac\gamma2\norm{f_F}^2,
\end{equation}
where
\begin{equation}\label{eq:weights}
 \epsilon=
 \begin{cases}
  \dfrac{\gamma a}{3a+(-\alpha)_++\gamma^2},&\text{ULD},\\[7pt]
  \dfrac{\gamma a}
        {100[\gamma^2+\min\{d(a+(-\alpha)_+),\,a+(-\alpha)_++\beta d\}]},
       &\text{BPS},\\[7pt]
  \dfrac{\gamma a}{100[a+\beta+(-\alpha)_++\gamma^2]},
       &\text{ZZ}.
 \end{cases}
\end{equation}
Moreover, the weight \eqref{eq:weights} satisfies
$0<\epsilon\le\min\{\gamma,\sqrt a/2\}$.
\end{proposition}

Subsection~\ref{sec:cross-terms} proves
Proposition~\ref{prop:common-dissipation} from the velocity dissipation
and sampler-dependent cross-term estimates.
Subsection~\ref{sec:resolvent} bounds $\mc D_a$ below using the LSI.
Subsection~\ref{sec:auxiliary-proofs} contains the sampler-specific
calculations.

\subsection{Velocity dissipation and cross terms}\label{sec:cross-terms}

The following velocity dissipation estimates are taken from
Fan, Li, and Lu~\cite{FanLiLu26}*{(25), (30)} for ULD and
Lu and Wang~\cite{LuWangPDMP}*{(30)--(32)} for BPS and ZZ.

\begin{lemma}[Velocity dissipation]
\label{lem:velocity-dissipation}
Let $f\in L^2_0(\mu)$ be a test function. For ULD,
\begin{equation}\label{eq:diss-uld}
 \ip{f}{\mc Lf}=-\gamma\norm{\nabla_vf}^2\le-\gamma\norm{f_F}^2,
\end{equation}
where the inequality is the Gaussian Poincar\'e inequality in~$v$.

For BPS,
\begin{equation}\label{eq:diss-bps}
 \begin{aligned}
 \ip{f}{\mc Lf}
 &=-\frac12\int(v\cdot\nabla U)_+
       [f(x,R_xv)-f(x,v)]^2\ud\mu
       -\gamma\norm{f_F}^2\\
 &\le-\gamma\norm{f_F}^2.
 \end{aligned}
\end{equation}

For ZZ,
\begin{equation}\label{eq:diss-zz}
 \begin{aligned}
 \ip{f}{\mc Lf}
 &=-\frac12\sum_{i=1}^d\int(v_i\partial_{x_i}U)_+
       [f(x,F_iv)-f(x,v)]^2\ud\mu
       -\gamma\norm{f_F}^2\\
 &\le-\gamma\norm{f_F}^2.
 \end{aligned}
\end{equation}
\end{lemma}

One term in the derivative of the correction in
\eqref{eq:modified-energy} is the pairing
$\ip{f}{\mc L(v\cdot\nabla_x\phi)}$, where
$\phi=(a-\mc L_{\mathrm o})^{-1}f_S$ is a function of~$x$ alone.
The velocity average of $\mc L(v\cdot\nabla_x\phi)$ pairs with~$f_S$
and yields the macroscopic dissipation~\eqref{eq:Q}. Its fluctuation
pairs with~$f_F$ and yields the cross term
$\ip{f_F}{(I-\Pi_v)\mc L(v\cdot\nabla_x\phi)}$.
The following estimates control this cross term by the macroscopic and velocity dissipations, so that it can be absorbed in the proof of Proposition~\ref{prop:common-dissipation}.

\begin{lemma}[Velocity-average identity and cross-term estimates]
\label{lem:model-estimates}
For $\phi\in C_c^\infty(\RR^d)$, each of the three generators satisfies
the velocity-average identity
\begin{equation}\label{eq:generator-average}
 \Pi_v \mc L(v\cdot\nabla_x\phi)=\mc L_{\mathrm o}\phi.
\end{equation}
For the remaining assertions, assume \eqref{eq:hessian-bounds},
with $\beta<\infty$ for ZZ.
Let $a>0$, let $f\in L^2_0(\mu)$ be a test function, and let
$\phi=(a-\mc L_{\mathrm o})^{-1}f_S$ be the associated resolvent function.
Then \eqref{eq:generator-average} holds for this~$\phi$, and we
have the following estimates:
\begin{enumerate}
\item For ULD,
\begin{equation}\label{eq:cross-uld}
 \left|\ip{f_F}{(I-\Pi_v)\mc L(v\cdot\nabla_x\phi)}\right|
 \le\sqrt{1+\frac{(-\alpha)_++\gamma^2}{a}}\,
       \sqrt{\mc D_a(f_S)}\norm{\nabla_vf}.
\end{equation}
\item For BPS,
\begin{equation}\label{eq:cross-bps}
 \begin{aligned}
 &\left|\ip{f_F}{(I-\Pi_v)\mc L(v\cdot\nabla_x\phi)}\right|\\
 &\quad\le\sqrt{66d+\frac{(30d+6)(-\alpha)_++3\gamma^2}{a}}\,
       \sqrt{\mc D_a(f_S)}\norm{f_F}.
 \end{aligned}
\end{equation}
If in addition $\beta<\infty$, then
\begin{equation}\label{eq:nonconvex-cross-bps-second}
 \begin{aligned}
 &\left|\ip{f_F}{(I-\Pi_v)\mc L(v\cdot\nabla_x\phi)}\right|\\
 &\quad\le\sqrt{66+\frac{66(-\alpha)_++30\beta d+3\gamma^2}{a}}\,
       \sqrt{\mc D_a(f_S)}\norm{f_F}.
 \end{aligned}
\end{equation}
\item For ZZ,
\begin{equation}\label{eq:cross-zz}
 \left|\ip{f_F}{(I-\Pi_v)\mc L(v\cdot\nabla_x\phi)}\right|
 \le\sqrt{66+\frac{66(-\alpha)_++30\beta+3\gamma^2}{a}}\,
       \sqrt{\mc D_a(f_S)}\norm{f_F}.
\end{equation}
\end{enumerate}
\end{lemma}
Lemma~\ref{lem:model-estimates} is proved for each sampler in
Subsections~\ref{sec:uld-estimates}--\ref{sec:zz-estimates}.

From now on, $\phi=(a-\mc L_{\mathrm o})^{-1}f_S$ denotes the
resolvent function of Lemma~\ref{lem:model-estimates}; it is defined for
every $f_S\in L^2_0(\mu_x)$. The sampler-specific proofs bound squared
norms of derivatives of~$\phi$ in terms of $\mc D_a(f_S)$. Pairing the resolvent equation
$f_S=(a-\mc L_{\mathrm o})\phi$ with $-\mc L_{\mathrm o}\phi$, and using
\eqref{eq:Q} and $\ip{\phi}{-\mc L_{\mathrm o}\phi}=\norm{\nabla_x\phi}^2$,
gives
\begin{equation}\label{eq:resolvent-energy}
 \mc D_a(f_S)
 =a\norm{\nabla_x\phi}^2+\norm{\mc L_{\mathrm o}\phi}^2.
\end{equation}
Both terms on the right are nonnegative, so
\begin{equation}\label{eq:resolvent-bounds}
 \norm{\mc L_{\mathrm o}\phi}^2\le \mc D_a(f_S),
 \qquad
 \norm{\nabla_x\phi}^2\le\frac1a\mc D_a(f_S).
\end{equation}
By \eqref{eq:Q}, Cauchy--Schwarz, and the bound on
$\norm{\mc L_{\mathrm o}\phi}^2$ in \eqref{eq:resolvent-bounds},
$\mc D_a(f_S)=\ip{f_S}{-\mc L_{\mathrm o}\phi}
\le\norm{f_S}\sqrt{\mc D_a(f_S)}$.
Hence $0\le\mc D_a(f_S)\le\norm{f_S}^2$, and $\mc D_a$ is a bounded
quadratic form on~$L^2_0(\mu_x)$.

\begin{proof}[Proof of Proposition~\ref{prop:common-dissipation}]
We first differentiate the modified energy along the forward
equation~\eqref{eq:density-evolution}. The calculation follows Fan, Li,
and Lu~\cite{FanLiLu26}*{(41)--(48)}, with the sign change explained in
Section~\ref{sec:modified-energy}. The forward equation gives
\begin{equation}\label{eq:uncorrected-derivative}
 \frac12\frac{\ud}{\ud t}\norm f^2
 =\ip{f}{\mc L^*f}.
\end{equation}
For test functions, $\ip{f}{\mc L^*f}=\ip{\mc Lf}{f}=\ip{f}{\mc Lf}$.
Since $C_c^\infty(\RR^{2d})$ is a core for $\mc L^*$, the estimates in
Lemmas~\ref{lem:velocity-dissipation} and~\ref{lem:model-estimates}
extend by approximation to $f\in D(\mc L^*)$, with $\ip{f}{\mc Lf}$
replaced by $\ip{f}{\mc L^*f}$. For the correction, the product rule
gives
\begin{equation}\label{eq:corrector-product-rule}
 \frac{\ud}{\ud t}\ip{\mc A_a f}{f}
 =\ip{\mc A_a\partial_tf}{f}
  +\ip{\mc A_a f}{\partial_tf}.
\end{equation}
Using $\Pi_v\mc A_a=\mc A_a$, the adjoint formula from
Lemma~\ref{lem:corrector-bounds}, and $\partial_tf=\mc L^*f$, we write the
first term as
\begin{align*}
 \ip{\mc A_a\partial_tf}{f}
 &=\ip{\mc A_a\partial_tf}{f_S}
 =\ip{\partial_tf}{\mc A_a^*f_S}\\
 &=\ip{\partial_tf}{v\cdot\nabla_x\phi}
 =\ip{f}{\mc L(v\cdot\nabla_x\phi)}.
\end{align*}
Since $\Pi_v f_S=f_S$, $\Pi_v f_F=0$, and $\Pi_v$ is self-adjoint, the
right-hand side splits as in the first line below. The second line
uses the velocity-average identity for~$\phi$ from
Lemma~\ref{lem:model-estimates}, and the third uses the definition
\eqref{eq:Q} of~$\mc D_a$:
\begin{align*}
 \ip{f}{\mc L(v\cdot\nabla_x\phi)}
 &=\ip{f_S}{\Pi_v\mc L(v\cdot\nabla_x\phi)}
   +\ip{f_F}{(I-\Pi_v)\mc L(v\cdot\nabla_x\phi)}\\
 &=\ip{f_S}{\mc L_{\mathrm o}\phi}
   +\ip{f_F}{(I-\Pi_v)\mc L(v\cdot\nabla_x\phi)}\\
 &=-\mc D_a(f_S)
   +\ip{f_F}{(I-\Pi_v)\mc L(v\cdot\nabla_x\phi)}.
\end{align*}
For the second term in \eqref{eq:corrector-product-rule}, the function
$\mc A_a f=\mc A_a f_F$ depends on~$x$ alone. Hence
$\mc L\mc A_af_F=(v\cdot\nabla_x)\mc A_af_F$, which has zero velocity
average, and
\begin{align*}
 \ip{\mc A_a f}{\partial_tf}
 &=\ip{\mc A_a f_F}{\mc L^* f}
 =\ip{\mc L\mc A_a f_F}{f}\\
 &=\ip{(v\cdot\nabla_x)\mc A_a f_F}{f}
 =\ip{(v\cdot\nabla_x)\mc A_a f_F}{f_F}.
\end{align*}
Together, these identities give
\begin{equation}\label{eq:corrector-derivative}
 \begin{aligned}
 \frac{\ud}{\ud t}\ip{\mc A_a f}{f}
 ={}&-\mc D_a(f_S)
   +\ip{f_F}{(I-\Pi_v)\mc L(v\cdot\nabla_x\phi)}\\
   &+\ip{f_F}{(v\cdot\nabla_x)\mc A_a f_F}.
 \end{aligned}
\end{equation}
Lemma~\ref{lem:corrector-bounds} controls the final term:
\begin{equation}\label{eq:last-corrector-term}
 \left|\ip{f_F}{(v\cdot\nabla_x)\mc A_a f_F}\right|
 \le\norm{f_F}\norm{(v\cdot\nabla_x)\mc A_a f_F}
 \le\norm{f_F}^2.
\end{equation}
Combining \eqref{eq:uncorrected-derivative},
\eqref{eq:corrector-derivative}, and \eqref{eq:last-corrector-term},
we obtain
\begin{equation}\label{eq:energy-derivative}
 \begin{aligned}
 \frac{\ud}{\ud t}\mf L_{a,\epsilon}(f)
 \le{}&\ip{f}{\mc L^*f}-\epsilon \mc D_a(f_S)\\
 &+\epsilon\ip{f_F}{(I-\Pi_v)\mc L(v\cdot\nabla_x\phi)}
   +\epsilon\norm{f_F}^2.
 \end{aligned}
\end{equation}

We now absorb the cross term. First,
any finite upper Hessian bound $\beta$ is positive: otherwise the
restriction of the coercive potential $U$ to every line would be concave,
which is impossible for a function tending to $+\infty$ in both directions.
Each weight in \eqref{eq:weights} consequently satisfies
\[
 0<\epsilon\le\frac{\gamma a}{a+\gamma^2}
 \le\min\{\gamma,\sqrt a/2\}.
\]
In particular, Lemma~\ref{lem:energy-comparison} applies.

Write $\mc D_v=\norm{\nabla_vf}^2$ for ULD and $\mc D_v=\norm{f_F}^2$ for BPS and ZZ, and let $K$ be the coefficient multiplying $\sqrt{\mc D_a(f_S)}\sqrt{\mc D_v}$ in the corresponding estimate of Lemma~\ref{lem:model-estimates}. For BPS, use the estimate corresponding to the smaller of the two terms inside the minimum in~\eqref{eq:weights}.
Lemmas~\ref{lem:velocity-dissipation} and~\ref{lem:model-estimates} give
$\ip{f}{\mc L^*f}\le-\gamma\mc D_v$ and $\norm{f_F}^2\le\mc D_v$,
together with the cross-term bound
$K\sqrt{\mc D_a(f_S)}\sqrt{\mc D_v}$.
Young's inequality in \eqref{eq:energy-derivative} yields
\begin{align*}
 \frac{\ud}{\ud t}\mf L_{a,\epsilon}(f)
 &\le-\epsilon \mc D_a(f_S)-\gamma\mc D_v
       +\epsilon K\sqrt{\mc D_a(f_S)}\sqrt{\mc D_v}+\epsilon\mc D_v\\
 &\le-\frac\epsilon2\mc D_a(f_S)
       -\left[\gamma-\epsilon\left(1+\frac{K^2}{2}\right)\right]\mc D_v.
\end{align*}
It remains to check $\epsilon(K^2+2)\le\gamma$.
The prefactors in Lemma~\ref{lem:model-estimates} satisfy
\[
 a(K^2+2)\le
 \begin{cases}
  3a+(-\alpha)_++\gamma^2,&\text{ULD},\\
  100[\gamma^2+d(a+(-\alpha)_+)],&\text{BPS, \eqref{eq:cross-bps}},\\
  100[\gamma^2+a+(-\alpha)_++\beta d],
    &\text{BPS, \eqref{eq:nonconvex-cross-bps-second}},\\
  100[\gamma^2+a+(-\alpha)_++\beta],&\text{ZZ}.
 \end{cases}
\]
For \eqref{eq:cross-bps}, use $66d+2\le68d$ and $30d+6\le36d$;
all other coefficients are already at most $100$. The corresponding
right-hand side equals $\gamma a/\epsilon$ by \eqref{eq:weights}.
Thus the remaining velocity coefficient is at least $\gamma/2$.
Since $\mc D_v\ge\norm{f_F}^2$, this proves \eqref{eq:common-dissipation}.
\end{proof}

\subsection{A nonlinear spatial resolvent estimate}
\label{sec:resolvent}

We now use the LSI to bound the macroscopic dissipation $\mc D_a$,
defined in \eqref{eq:Q}, from below. The Poincar\'e inequality implied
by~\eqref{eq:lsi} and spectral calculus give the linear estimate
\begin{equation}\label{eq:linear-gap}
 \mc D_a(f_S)\ge\frac{\CLSI}{a+\CLSI}\norm{f_S}^2.
\end{equation}
Together with Proposition~\ref{prop:common-dissipation}, it yields
exponential decay of~$\chi^2$, and hence mixing times linear in~$R$. When
$1+f_S$ is a probability density, the LSI gives a nonlinear estimate
instead. Its numerator has the structure needed for the R\'enyi-2
divergence: if $\norm f^2$ decreased at a rate at least
$c(1+\norm f^2)\log(1+\norm f^2)$, then $\log(1+\norm f^2)$ would decay
exponentially at rate~$c$.

\begin{proposition}[Nonlinear spatial resolvent estimate]
\label{prop:nonlinear-resolvent}
Assume \eqref{eq:lsi}. For every $a>0$ and every
$f_S\in L^2_0(\mu_x)$ satisfying $f_S\ge-1$,
\begin{equation}\label{eq:nonlinear-resolvent-general}
 \mc D_a(f_S)
 \ge\frac{\CLSI(1+\norm{f_S}^2)\log(1+\norm{f_S}^2)}
 {2a+\CLSI(1+\log(1+\norm{f_S}^2))}.
\end{equation}
\end{proposition}

\begin{proof}
We apply the LSI to a density obtained by smoothing $1+f_S$ with the normalized square root of the resolvent. The resulting density has squared gradient norm $a\mc D_a(f_S)$. Define
\begin{equation}\label{eq:square-root-mixture}
 \psi=\sqrt a(a-\mc L_{\mathrm o})^{-1/2}(1+f_S)
 =\sqrt{\frac a\pi}\int_0^\infty
       t^{-1/2}e^{-at}e^{t\mc L_{\mathrm o}}(1+f_S)\ud t.
\end{equation}
The second equality follows from spectral calculus, since the spectrum
of $a-\mc L_{\mathrm o}$ lies in $[a,\infty)$, and the identity
\[
 \lambda^{-1/2}=\frac1{\sqrt\pi}\int_0^\infty t^{-1/2}e^{-\lambda t}\ud t,
 \qquad\lambda>0.
\]
Taking $\lambda=a$ shows that the weight
$\sqrt{a/\pi}\,t^{-1/2}e^{-at}$ has integral one. The semigroup
$e^{t\mc L_{\mathrm o}}$ preserves probability densities and contracts
$L^2(\mu_x)$. Thus the integral converges in $L^2(\mu_x)$ and defines a
probability density~$\psi$.

Since $\mc L_{\mathrm o}1=0$, the definition gives
$\psi=1+\sqrt a\,(a-\mc L_{\mathrm o})^{-1/2}f_S$, and the second term
has mean zero because $f_S$ does. The operator
$(-\mc L_{\mathrm o})^{1/2}(a-\mc L_{\mathrm o})^{-1/2}$ is bounded, so
$\psi$ lies in the domain of the form $h\mapsto\norm{\nabla_xh}^2$.
Using \eqref{eq:Q} and
$a(a-\mc L_{\mathrm o})^{-1}=I-(-\mc L_{\mathrm o})(a-\mc L_{\mathrm o})^{-1}$,
we obtain
\begin{equation}\label{eq:square-root-identity}
 \begin{aligned}
 \norm{\nabla_x\psi}^2
 &=a\ip{f_S}{(-\mc L_{\mathrm o})(a-\mc L_{\mathrm o})^{-1}f_S}
 =a\mc D_a(f_S),\\
 \norm\psi^2
 &=1+a\ip{f_S}{(a-\mc L_{\mathrm o})^{-1}f_S}
 =1+\norm{f_S}^2-\mc D_a(f_S).
 \end{aligned}
\end{equation}
We next bound $\Ent_{\mu_x}(\psi^2)$ from both sides. Since $\psi$ is a
probability density, Jensen's inequality for $u\mapsto u\log u$ under
the probability measure $\psi\ud\mu_x$ gives
\[
 \int\psi^2\log\psi^2\ud\mu_x
 =2\int\psi\,(\psi\log\psi)\ud\mu_x
 \ge2\norm\psi^2\log\norm\psi^2.
\]
Subtracting $\norm\psi^2\log\norm\psi^2$ gives the lower bound below.
The upper bound follows from the LSI and the identity for
$\norm{\nabla_x\psi}^2$ in \eqref{eq:square-root-identity}:
\begin{equation}\label{eq:square-root-entropy}
 \norm\psi^2\log\norm\psi^2
 \le\Ent_{\mu_x}(\psi^2)
 =\int\psi^2\log\psi^2\ud\mu_x-\norm\psi^2\log\norm\psi^2
 \le\frac{2a}{\CLSI}\mc D_a(f_S).
\end{equation}

The convexity of $u\log u$, the identity for $\norm\psi^2$ in
\eqref{eq:square-root-identity}, and \eqref{eq:square-root-entropy} give
\begin{align*}
 (1+\norm{f_S}^2)\log(1+\norm{f_S}^2)
 &\le\norm\psi^2\log\norm\psi^2
   +\bigl[1+\log(1+\norm{f_S}^2)\bigr]\mc D_a(f_S)\\
 &\le\left[\frac{2a}{\CLSI}+1+\log(1+\norm{f_S}^2)\right]\mc D_a(f_S).
\end{align*}
Rearranging proves \eqref{eq:nonlinear-resolvent-general}.
\end{proof}

\subsection{Proofs of the sampler-specific estimates}\label{sec:auxiliary-proofs}

We now prove Lemma~\ref{lem:model-estimates} for each sampler. All
calculations are initially on test functions and extend to the resolvent
functions by standard approximation.

\subsubsection{The Bochner identity, a weighted estimate, and Gaussian moments}
\label{sec:common-estimates}

The integrated Bochner identity~\cite{FanLiLu26}*{(31)} reads
\begin{equation}\label{eq:bochner}
 \norm{\mc L_{\mathrm o}\phi}^2
 =\norm{\nabla_x^2\phi}^2
   +\int\nabla_x^2U[\nabla_x\phi,\nabla_x\phi]\ud\mu_x.
\end{equation}
The lower Hessian bound therefore gives
\begin{equation}\label{eq:nonconvex-bochner}
 \norm{\nabla_x^2\phi}^2
 \le\norm{\mc L_{\mathrm o}\phi}^2
       +(-\alpha)_+\norm{\nabla_x\phi}^2.
\end{equation}
The additional term records the contribution of negative curvature. We retain it in each sampler-specific estimate below; it vanishes when~$\alpha\ge0$. Together with~\eqref{eq:resolvent-bounds}, this estimate applies for every~$a>0$, without a restriction~$a>(-\alpha)_+$.

When $\beta<\infty$, the two PDMPs will use the same weighted
coordinate estimate. For $u\in C_c^\infty(\RR^d)$ and each $i$,
\begin{equation}\label{eq:coordinate-weight}
 \norm{(\partial_{x_i}U)u}^2
 \le2\beta\norm u^2+4\norm{\partial_{x_i}u}^2.
\end{equation}
Integration by parts and Young's inequality give
\begin{align*}
 \norm{(\partial_{x_i}U)u}^2
 &=\int(\partial_{x_ix_i}U)u^2\ud\mu_x
   +2\int(\partial_{x_i}U)u\partial_{x_i}u\ud\mu_x\\
 &\le\beta\norm u^2
   +\frac12\norm{(\partial_{x_i}U)u}^2
   +2\norm{\partial_{x_i}u}^2.
\end{align*}
Absorbing the repeated term proves \eqref{eq:coordinate-weight}.

We use the elementary Gaussian moment identities
\[
 \int v_iv_j\ud\mu_v=\delta_{ij},\qquad
 \int v_iv_jv_kv_l\ud\mu_v
 =\delta_{ij}\delta_{kl}+\delta_{ik}\delta_{jl}
  +\delta_{il}\delta_{jk}.
\]
Since $\nabla_x^2\phi$ is symmetric, for each $x$ this implies
\begin{align*}
 &\int(v^\top\nabla_x^2\phi\,v-\Delta_x\phi)^2\ud\mu_v\\
 &\quad=\sum_{i,j,k,l=1}^d
  (\partial_{x_ix_j}\phi)(\partial_{x_kx_l}\phi)
  \int(v_iv_j-\delta_{ij})(v_kv_l-\delta_{kl})\ud\mu_v\\
 &\quad=2\sum_{i,j=1}^d|\partial_{x_ix_j}\phi|^2
 =2|\nabla_x^2\phi|_{\mathrm{HS}}^2.
\end{align*}
Likewise, $\int(v\cdot\nabla_x\phi)^2\ud\mu_v=|\nabla_x\phi|^2$.
For unit $n$, the identities
$\int(v\cdot n)_+^2\ud\mu_v=1/2$ and
$\int(v\cdot n)_+^4\ud\mu_v=3/2$ give
\[
 \int\bigl(2(v\cdot n)_+^2-1\bigr)^2\ud\mu_v
 =4\cdot\frac32-4\cdot\frac12+1=5.
\]

\subsubsection{Estimates for underdamped Langevin dynamics}
\label{sec:uld-estimates}

\begin{proof}[Proof of Lemma~\ref{lem:model-estimates}, ULD part]
Applying each part of \eqref{eq:uld} to $v\cdot\nabla_x\phi$, which is
linear in $v$, gives
\begin{align*}
 (v\cdot\nabla_x)(v\cdot\nabla_x\phi)
   &=v^\top\nabla_x^2\phi\,v,\\
 (-\nabla U\cdot\nabla_v)(v\cdot\nabla_x\phi)
   &=-\nabla U\cdot\nabla_x\phi,\\
 (\Delta_v-v\cdot\nabla_v)(v\cdot\nabla_x\phi)
   &=-v\cdot\nabla_x\phi.
\end{align*}
Since $\int v_i\ud\mu_v=0$ and
$\int v_iv_j\ud\mu_v=\delta_{ij}$, averaging the sum proves
\eqref{eq:generator-average}. Subtracting this average gives
\begin{equation}\label{eq:uld-fluctuation}
 (I-\Pi_v)\mc L(v\cdot\nabla_x\phi)
 =v^\top\nabla_x^2\phi\,v-\Delta_x\phi
  -\gamma v\cdot\nabla_x\phi.
\end{equation}
In particular, the force term disappears under $I-\Pi_v$.

To estimate the pairing with $f_F$, integrate once in velocity. For
each $x$, Gaussian integration by parts gives
\begin{align*}
 \int f_F(v_iv_j-\delta_{ij})\ud\mu_v
 &=\int\partial_{v_i}(v_jf_F)\ud\mu_v
   -\delta_{ij}\int f_F\ud\mu_v
 =\int v_j\partial_{v_i}f\ud\mu_v,\\
 \int f_Fv_i\ud\mu_v&=\int\partial_{v_i}f\ud\mu_v.
\end{align*}
Multiplying the first identity by $\partial_{x_ix_j}\phi$ and the
second by $-\gamma\partial_{x_i}\phi$, summing, and integrating in $x$
shows that
\begin{equation}\label{eq:uld-velocity-ibp}
 \ip{f_F}{(I-\Pi_v)\mc L(v\cdot\nabla_x\phi)}
 =\ip{\nabla_vf}{(\nabla_x^2\phi)v-\gamma\nabla_x\phi}.
\end{equation}
The Gaussian first and second moments also imply
\begin{align*}
 &\norm{(\nabla_x^2\phi)v-\gamma\nabla_x\phi}^2\\
 &\quad=\int\sum_{i,j,k=1}^d
       (\partial_{x_ix_j}\phi)(\partial_{x_ix_k}\phi)v_jv_k\ud\mu
    +\gamma^2\norm{\nabla_x\phi}^2\\
 &\qquad\quad
    -2\gamma\int\sum_{i,j=1}^d
       (\partial_{x_ix_j}\phi)(\partial_{x_i}\phi)v_j\ud\mu\\
 &\quad=\norm{\nabla_x^2\phi}^2
      +\gamma^2\norm{\nabla_x\phi}^2.
\end{align*}
The Hessian bound~\eqref{eq:nonconvex-bochner}, the resolvent estimates~\eqref{eq:resolvent-bounds}, and Cauchy--Schwarz in
\eqref{eq:uld-velocity-ibp} give
\begin{align*}
 \left|\ip{f_F}{(I-\Pi_v)\mc L(v\cdot\nabla_x\phi)}\right|
 &\le\norm{\nabla_vf}
  \sqrt{\norm{\nabla_x^2\phi}^2+\gamma^2\norm{\nabla_x\phi}^2}
  \\
 &\le\sqrt{1+\frac{(-\alpha)_++\gamma^2}{a}}\,
       \sqrt{\mc D_a(f_S)}\norm{\nabla_vf}.
\end{align*}
This proves \eqref{eq:cross-uld}.
\end{proof}

\subsubsection{Estimates for the bouncy particle sampler}
\label{sec:bps-estimates}

\begin{proof}[Proof of Lemma~\ref{lem:model-estimates}, BPS part]
The reflection formula gives
\[
 (R_xv)\cdot\nabla_x\phi-v\cdot\nabla_x\phi
 =-2(v\cdot n)(n\cdot\nabla_x\phi).
\]
Since $\nabla U=|\nabla U|n$, multiplication by the bounce rate yields
\[
 (v\cdot\nabla U)_+
   [(R_xv)\cdot\nabla_x\phi-v\cdot\nabla_x\phi]
 =-2(v\cdot n)_+^2\nabla U\cdot\nabla_x\phi.
\]
The transport and refreshment terms are computed directly, giving
\begin{equation}\label{eq:bps-linear-test}
 \mc L(v\cdot\nabla_x\phi)
 =v^\top\nabla_x^2\phi\,v
   -2(v\cdot n)_+^2\nabla U\cdot\nabla_x\phi
   -\gamma v\cdot\nabla_x\phi.
\end{equation}
Since $\int(v\cdot n)_+^2\ud\mu_v=1/2$ for unit $n$
(Section~\ref{sec:common-estimates}), the velocity average of
\eqref{eq:bps-linear-test} is
$\Delta_x\phi-\nabla U\cdot\nabla_x\phi$. This proves \eqref{eq:generator-average} for BPS and yields
\begin{equation}\label{eq:bps-centered-test}
 \begin{aligned}
 (I-\Pi_v)\mc L(v\cdot\nabla_x\phi)
 ={}&v^\top\nabla_x^2\phi\,v-\Delta_x\phi\\
 &-\bigl(2(v\cdot n)_+^2-1\bigr)\nabla U\cdot\nabla_x\phi
  -\gamma v\cdot\nabla_x\phi.
 \end{aligned}
\end{equation}

The Gaussian moments in Section~\ref{sec:common-estimates} control
the three terms. Using $(r+s+t)^2\le3(r^2+s^2+t^2)$ in
\eqref{eq:bps-centered-test} and integrating proves
\begin{equation}\label{eq:bps-fluctuation}
 \norm{(I-\Pi_v)\mc L(v\cdot\nabla_x\phi)}^2
 \le6\norm{\nabla_x^2\phi}^2
   +15\norm{\nabla U\cdot\nabla_x\phi}^2
   +3\gamma^2\norm{\nabla_x\phi}^2.
\end{equation}

We bound the force term $\norm{\nabla U\cdot\nabla_x\phi}^2$ in two
ways, leading to \eqref{eq:cross-bps} and
\eqref{eq:nonconvex-cross-bps-second}. For \eqref{eq:cross-bps}, we use
only the lower Hessian bound~\eqref{eq:nonconvex-bochner}. Pointwise,
\[
 |\Delta_x\phi|^2
 =\left|\sum_{i=1}^d\partial_{x_ix_i}\phi\right|^2
 \le d\sum_{i=1}^d|\partial_{x_ix_i}\phi|^2
 \le d|\nabla_x^2\phi|_{\mathrm{HS}}^2.
\]
Using $\nabla U\cdot\nabla_x\phi
=\Delta_x\phi-\mc L_{\mathrm o}\phi$, we obtain
\begin{align}
 \norm{\nabla U\cdot\nabla_x\phi}^2
 &\le2\norm{\Delta_x\phi}^2+2\norm{\mc L_{\mathrm o}\phi}^2
 \nonumber\\
 &\le2d\norm{\nabla_x^2\phi}^2+2\norm{\mc L_{\mathrm o}\phi}^2\nonumber\\
 &\le2(d+1)\norm{\mc L_{\mathrm o}\phi}^2
       +2d(-\alpha)_+\norm{\nabla_x\phi}^2.
 \label{eq:force-one}
\end{align}
Substituting in \eqref{eq:bps-fluctuation}, and then using
\eqref{eq:resolvent-bounds}, gives
\begin{align}
 \norm{(I-\Pi_v)\mc L(v\cdot\nabla_x\phi)}^2
 &\le(30d+36)\norm{\mc L_{\mathrm o}\phi}^2\nonumber\\
 &\quad+[(30d+6)(-\alpha)_++3\gamma^2]\norm{\nabla_x\phi}^2\nonumber\\
 &\le\left(66d+\frac{(30d+6)(-\alpha)_++3\gamma^2}{a}\right)
           \mc D_a(f_S),
 \label{eq:bps-fluctuation-one}
\end{align}
using $d\ge1$ in the last inequality.

For \eqref{eq:nonconvex-cross-bps-second}, assume also $\beta<\infty$.
The first inequality below is the pointwise Cauchy--Schwarz bound
$|\nabla U\cdot\nabla_x\phi|^2\le|\nabla U|^2|\nabla_x\phi|^2$.
Then apply \eqref{eq:coordinate-weight} with $u=\partial_{x_j}\phi$, sum
over $i,j$, and use \eqref{eq:nonconvex-bochner}:
\begin{align}
 \norm{\nabla U\cdot\nabla_x\phi}^2
 &\le\sum_{i,j=1}^d
       \norm{(\partial_{x_i}U)(\partial_{x_j}\phi)}^2\nonumber\\
 &\le2\beta d\norm{\nabla_x\phi}^2+4\norm{\nabla_x^2\phi}^2\nonumber\\
 &\le4\norm{\mc L_{\mathrm o}\phi}^2
       +[2\beta d+4(-\alpha)_+]\norm{\nabla_x\phi}^2.
 \label{eq:force-two}
\end{align}
This time \eqref{eq:bps-fluctuation} and
\eqref{eq:resolvent-bounds} give
\begin{align}
 \norm{(I-\Pi_v)\mc L(v\cdot\nabla_x\phi)}^2
 &\le66\norm{\mc L_{\mathrm o}\phi}^2\nonumber\\
 &\quad+[66(-\alpha)_++30\beta d+3\gamma^2]\norm{\nabla_x\phi}^2\nonumber\\
 &\le\left(66+\frac{66(-\alpha)_++30\beta d+3\gamma^2}{a}\right)
       \mc D_a(f_S).
 \label{eq:bps-fluctuation-two}
\end{align}
Thus Cauchy--Schwarz gives the cross-term estimates~\eqref{eq:cross-bps} and~\eqref{eq:nonconvex-cross-bps-second}.
\end{proof}

\subsubsection{Estimates for the Zigzag sampler}
\label{sec:zz-estimates}

\begin{proof}[Proof of Lemma~\ref{lem:model-estimates}, ZZ part]
For a coordinate flip,
\[
 (F_iv)\cdot\nabla_x\phi-v\cdot\nabla_x\phi
 =-2v_i\partial_{x_i}\phi.
\]
Consequently,
\begin{equation}\label{eq:zz-linear-test}
 \mc L(v\cdot\nabla_x\phi)
 =v^\top\nabla_x^2\phi\,v
  -2\sum_{i=1}^d(v_i\partial_{x_i}U)_+v_i\partial_{x_i}\phi
  -\gamma v\cdot\nabla_x\phi.
\end{equation}
Evenness of the one-dimensional Gaussian gives, for each fixed $x$,
\begin{align*}
 \int2v_i(v_i\partial_{x_i}U)_+\ud\mu_v
 &=\partial_{x_i}U,\\
 \int[2v_i(v_i\partial_{x_i}U)_+]^2\ud\mu_v
 &=6(\partial_{x_i}U)^2.
\end{align*}
Both identities also hold when $\partial_{x_i}U=0$.
The velocity average of \eqref{eq:zz-linear-test} is therefore
$\Delta_x\phi-\nabla U\cdot\nabla_x\phi$, proving
\eqref{eq:generator-average} for ZZ. Subtracting this average yields
\begin{equation}\label{eq:zz-centered-test}
 \begin{aligned}
 (I-\Pi_v)\mc L(v\cdot\nabla_x\phi)
 ={}&v^\top\nabla_x^2\phi\,v-\Delta_x\phi\\
 &-\sum_{i=1}^d
  [2v_i(v_i\partial_{x_i}U)_+-\partial_{x_i}U]
       \partial_{x_i}\phi
  -\gamma v\cdot\nabla_x\phi.
 \end{aligned}
\end{equation}
The summands in the second term have mean zero under~$\mu_v$ and depend
on independent Gaussian coordinates. Thus the two moment identities above imply
\begin{align*}
 &\int\left\{\sum_{i=1}^d
  [2v_i(v_i\partial_{x_i}U)_+-\partial_{x_i}U]
       \partial_{x_i}\phi\right\}^2\ud\mu_v\\
 &\qquad=5\sum_{i=1}^d
       (\partial_{x_i}U)^2(\partial_{x_i}\phi)^2.
\end{align*}
Using the quadratic Gaussian moment in
Section~\ref{sec:common-estimates} and
$(r+s+t)^2\le3(r^2+s^2+t^2)$ in
\eqref{eq:zz-centered-test}, we obtain
\begin{equation}\label{eq:zz-fluctuation}
 \begin{aligned}
 \norm{(I-\Pi_v)\mc L(v\cdot\nabla_x\phi)}^2
 \le{}&6\norm{\nabla_x^2\phi}^2
       +15\sum_{i=1}^d
        \norm{(\partial_{x_i}U)(\partial_{x_i}\phi)}^2\\
 &+3\gamma^2\norm{\nabla_x\phi}^2.
 \end{aligned}
\end{equation}

Here \eqref{eq:coordinate-weight} is needed only with
$u=\partial_{x_i}\phi$. Summing over $i$, rather than over $i,j$ as
for BPS, avoids the dimension factor:
\begin{align}
 \sum_{i=1}^d\norm{(\partial_{x_i}U)(\partial_{x_i}\phi)}^2
 &\le2\beta\norm{\nabla_x\phi}^2
       +4\sum_{i=1}^d\norm{\partial_{x_ix_i}\phi}^2\nonumber\\
 &\le4\norm{\mc L_{\mathrm o}\phi}^2
       +[2\beta+4(-\alpha)_+]\norm{\nabla_x\phi}^2.
 \label{eq:zz-coordinate-force}
\end{align}
The last step uses
$\sum_i\norm{\partial_{x_ix_i}\phi}^2\le\norm{\nabla_x^2\phi}^2$
and \eqref{eq:nonconvex-bochner}.
Combining \eqref{eq:zz-fluctuation},
\eqref{eq:zz-coordinate-force}, and \eqref{eq:resolvent-bounds},
we conclude that
\begin{align}
 \norm{(I-\Pi_v)\mc L(v\cdot\nabla_x\phi)}^2
 &\le66\norm{\mc L_{\mathrm o}\phi}^2\nonumber\\
 &\quad+[66(-\alpha)_++30\beta+3\gamma^2]\norm{\nabla_x\phi}^2\nonumber\\
 &\le\left(66+\frac{66(-\alpha)_++30\beta+3\gamma^2}{a}\right)
       \mc D_a(f_S).
 \label{eq:zz-fluctuation-bound}
\end{align}
Hence Cauchy--Schwarz gives \eqref{eq:cross-zz}.
\end{proof}

\section{Convergence and mixing-time bounds}\label{sec:closure}

We prove the decay estimate for the R\'enyi-2 divergence and the resulting
mixing-time bounds.

\begin{proof}[Proof of Theorem~\ref{thm:main}]
We apply
Proposition~\ref{prop:nonlinear-resolvent} to the spatial component
$f_{S,t}$, which lies in $L^2_0(\mu_x)$ and satisfies $f_{S,t}\ge-1$
along each dynamics. Orthogonal projection and the
$L^2(\mu)$-contractivity of $\mc P_t^*$ give
\[
 0\le\log(1+\norm{f_S}^2)
 \le\log(1+\norm{f_t}^2)
 \le\log(1+\norm{f_0}^2)\le R.
\]
Choose $a=\CLSI(1+R)$ to control
the denominator in Proposition~\ref{prop:nonlinear-resolvent}
uniformly in time:
\[
 2a+\CLSI(1+\log(1+\norm{f_S}^2))
 \le2a+\CLSI(1+R)=3a.
\]
Thus, at every time,
\begin{equation}\label{eq:nonlinear-resolvent}
 \mc D_a(f_S)
 \ge\frac{\CLSI}{3a}(1+\norm{f_S}^2)\log(1+\norm{f_S}^2).
\end{equation}

We next combine the macroscopic and velocity dissipation into
\begin{align}\label{eq:full-density-dissipation}
 \frac{\ud}{\ud t}\mf L_{a,\epsilon}(f_t)
 &\le-\frac{\CLSI\epsilon}{6a}
       (1+\norm{f_S}^2)\log(1+\norm{f_S}^2)
       -\frac\gamma2\norm{f_F}^2\nonumber\\
 &\le-\frac{\CLSI\epsilon}{6a}
       \left[(1+\norm{f_S}^2)\log(1+\norm{f_S}^2)
              +(1+R)\norm{f_F}^2\right]\nonumber\\
 &\le-\frac{\CLSI\epsilon}{6a}
       (1+\norm{f_t}^2)\log(1+\norm{f_t}^2).
\end{align}
The first inequality follows from
Proposition~\ref{prop:common-dissipation} and
\eqref{eq:nonlinear-resolvent}. The second holds because
$a=\CLSI(1+R)$ and $\epsilon\le\gamma$ give
\[
 \frac{\CLSI\epsilon}{6a}(1+R)
 =\frac{\epsilon}{6}\le\frac\gamma2.
\]
For the last inequality, orthogonality gives
$\norm{f}^2-\norm{f_S}^2=\norm{f_F}^2$.
Both $1+\norm{f}^2$ and $1+\norm{f_S}^2$ are at least one, and their
logarithms are at most $R$, so
\begin{align*}
 &(1+\norm{f}^2)\log(1+\norm{f}^2)
   -(1+\norm{f_S}^2)\log(1+\norm{f_S}^2)\\
 &\qquad=\int_{1+\norm{f_S}^2}^{1+\norm{f}^2}
                   (1+\log s)\ud s\\
 &\qquad\le(1+R)
          \bigl(\norm{f}^2-\norm{f_S}^2\bigr)
 =(1+R)\norm{f_F}^2.
\end{align*}

It remains to turn \eqref{eq:full-density-dissipation} into
decay of $\Ren_2(\varrho_t\mid\mu)=\log(1+\norm{f_t}^2)$. The
right-hand side of \eqref{eq:full-density-dissipation} equals
$-\frac{\CLSI\epsilon}{6a}(1+\norm{f_t}^2)\,\Ren_2(\varrho_t\mid\mu)$,
but the left-hand side is the derivative of the modified energy rather
than of $\norm{f_t}^2$. We therefore replace $\norm{f_t}^2$ in
$\Ren_2$ by a multiple of the modified energy and set
\[
 H(t)=\log\Bigl(1+\frac43\mf L_{a,\epsilon}(f_t)\Bigr).
\]
The factor $4/3$ is chosen so that Lemma~\ref{lem:energy-comparison}
gives
\[
 \frac13\norm{f_t}^2\le\frac43\mf L_{a,\epsilon}(f_t)\le\norm{f_t}^2.
\]
Since $\log(1+x/3)\ge\frac13\log(1+x)$ for $x\ge0$ by concavity, it
follows that
\begin{equation}\label{eq:H-comparison}
 \frac13\Ren_2(\varrho_t\mid\mu)\le H(t)
 \le\Ren_2(\varrho_t\mid\mu).
\end{equation}

By \eqref{eq:full-density-dissipation},
\[
 H'(t)
 =\frac{\frac43\frac{\ud}{\ud t}\mf L_{a,\epsilon}(f_t)}
       {1+\frac43\mf L_{a,\epsilon}(f_t)}
 \le-\frac{2\CLSI\epsilon}{9a}\,
   \frac{1+\norm{f_t}^2}{1+\frac43\mf L_{a,\epsilon}(f_t)}\,
   \Ren_2(\varrho_t\mid\mu).
\]
The fraction $(1+\norm{f_t}^2)/(1+\frac43\mf L_{a,\epsilon}(f_t))$ is at
least one by the upper bound on $\frac43\mf L_{a,\epsilon}(f_t)$, and
$\Ren_2(\varrho_t\mid\mu)\ge H(t)$ by \eqref{eq:H-comparison}. Hence
\begin{equation}\label{eq:modified-log-decay}
 H'(t)\le-\frac{2\CLSI\epsilon}{9a}H(t).
\end{equation}
Gr\"onwall's inequality for \eqref{eq:modified-log-decay}, together with
\eqref{eq:H-comparison}, gives
\[
 \Ren_2(\varrho_t\mid\mu)
 \le3H(t)
 \le3e^{-2\CLSI\epsilon t/(9a)}H(0)
 \le3e^{-2\CLSI\epsilon t/(9a)}\Ren_2(\varrho_0\mid\mu).
\]
This is \eqref{eq:decay} with $\Lambda=2\CLSI\epsilon/(9a)$.
Substituting \eqref{eq:weights} and $a=\CLSI(1+R)$ gives
\eqref{eq:rate-uld}--\eqref{eq:rate-zz}.

This proves the theorem when $f_0\in D(\mc L^*)$, as required in
Proposition~\ref{prop:common-dissipation}. Standard approximation and
the $L^2(\mu)$-contractivity of $\mc P_t^*$ extend it to all initial laws
with finite R\'enyi-2 divergence.
\end{proof}

The choice of the shift drives the improved dependence on~$R$. In the
proof above, the rate is proportional to~$\CLSI\epsilon$ divided by the
denominator in Proposition~\ref{prop:nonlinear-resolvent}, which is at
most $2a+\CLSI(1+R)$ along the evolution. The shifted corrector has norm
at most~$1/(2\sqrt a)$, so the modified energy remains comparable to the
quadratic energy for every~$\epsilon\le\sqrt a$
(Lemma~\ref{lem:energy-comparison}). For fixed curvature and dimension,
\eqref{eq:weights} shows that a friction or refreshment rate~$\gamma$ of
order~$\sqrt a$ makes a weight of order~$\sqrt a$ admissible. The rate is
then of order~$\CLSI\sqrt a/(a+\CLSI R)$ for~$R\ge1$. In the linear
estimate~\eqref{eq:linear-gap}, the factor $\CLSI/(a+\CLSI)$ decreases as
$a$ grows, and it is best to take
$a=0$~\cite{FanLiLu26}*{Corollary~2}. With~$a$ of order~$\CLSI$, the
rate above is of order~$\sqrt{\CLSI}/R$. Taking~$a=\CLSI(1+R)$
raises it to order~$\sqrt{\CLSI/R}$, while the denominator
becomes~$3\CLSI(1+R)$, only a factor three above its zero-shift
value~$\CLSI(1+R)$.

\begin{proof}[Proof of Corollary~\ref{cor:mixing}]
We first prove the R\'enyi-2 bound, with $\gamma$ given
by~\eqref{eq:time-optimized-parameters}. We apply Theorem~\ref{thm:main} in stages with successively halved
divergence bounds, keeping $\gamma$ fixed. For a divergence bound $r$,
write the rate in Theorem~\ref{thm:main} as
\[
 \Lambda(r)=\frac{\CLSI\gamma}{c[\gamma^2+B(r)]},
\]
where $c=9/2$ for ULD and $c=450$ for BPS and ZZ.
Here $B(r)$ is the part of the corresponding bracket in
\eqref{eq:rate-uld}--\eqref{eq:rate-zz} independent of $\gamma^2$.
It is an affine function of $r$ with nonnegative coefficients for ULD
and ZZ, and the minimum of two such functions for BPS.
The choice \eqref{eq:time-optimized-parameters} is precisely
$\gamma^2=B(R/\mathfrak{L})$.

Set
\[
 N=\left\lceil\log_2\frac{R}{\eta}\right\rceil,
 \qquad R_j=2^{-j}R\quad(0\le j\le N),
 \qquad t_j=\frac{\log6}{\Lambda(R_j)}\quad(0\le j<N).
\]
If the divergence at the start of stage $j$ is at most $R_j$,
Theorem~\ref{thm:main} bounds it at the end by
\[
 3e^{-\Lambda(R_j)t_j}R_j
 =R_{j+1}.
\]
The theorem applies at each endpoint, including when $R_j<1$.
Induction therefore gives divergence at most $R_N\le\eta$
at time $\sum_{j<N}t_j$, and hence at all later times, since
$\norm{f_t}$ is nonincreasing.

Since $R/\eta\ge2$, the elementary bounds
\[
 \frac{\mathfrak{L}}{\log6}\le N\le\frac{\mathfrak{L}}{\log2},
 \qquad \sum_{j=0}^{N-1}R_j\le2R
\]
imply, for every affine function $b(r)=b_0+b_1r$ with $b_0,b_1\ge0$,
\[
 \sum_{j=0}^{N-1}b(R_j)
 \le Nb_0+2b_1R
 \le4N b(R/\mathfrak{L}).
\]
For BPS, apply this bound to each finite branch and take their minimum.
In all three cases,
\[
 \sum_{j=0}^{N-1}B(R_j)
 \le4N B(R/\mathfrak{L})=4N\gamma^2.
\]
It follows that
\begin{align*}
 \sum_{j=0}^{N-1}t_j
 &\le\frac{c\log6}{\CLSI\gamma}
       \left[N\gamma^2+\sum_{j=0}^{N-1}B(R_j)\right]\\
 &\le\frac{5c\log6}{\log2}\frac{\mathfrak{L}\gamma}{\CLSI}.
\end{align*}
The prefactor is below $60$ for ULD and $6000$ for BPS and ZZ.
Substituting \eqref{eq:time-optimized-parameters} proves the R\'enyi-2
bound.

For the total variation bound, $\eta=2\delta^2$ lies in $(0,1/2]$ for
$0<\delta\le1/2$, and~\eqref{eq:tv-renyi} gives
$\TV(\varrho_t,\mu)\le\delta$ whenever
$\Ren_2(\varrho_t\mid\mu)\le2\delta^2$.
\end{proof}

\appendix

\section{Sharpness for the R\'enyi-2 divergence}\label{app:initial-divergence-sharpness}

This appendix proves Theorem~\ref{thm:initial-divergence-sharpness}, which
shows that the order~$\sqrt{R\log R}$ in Corollary~\ref{cor:mixing} cannot
be improved for the R\'enyi-2 divergence of BPS and ZZ.

\begin{proof}[Proof of Theorem~\ref{thm:initial-divergence-sharpness}]
We construct one initial law, independent of~$\gamma$, with two
obstructions to the decay of the R\'enyi-2 divergence. Fix $\gamma>0$, and
suppose that $\Ren_2(\varrho_t\mid\mu)\le1/2$ at some time~$t\ge0$. Slow
refreshment retains a concentrated speed distribution, which gives a
lower bound on~$t$ of order~$R/\gamma$. Frequent refreshment delays the
return from a distant initial position, which gives a lower bound of
order~$\gamma\log R/\alpha$. The larger of the two bounds is at least
their geometric mean, which does not depend on~$\gamma$.

\emph{Initial law.} For~$U(x)=\alpha x^2/2$, the BPS reflection and the
ZZ flip both send~$v$ to~$-v$ at rate~$(\alpha xv)_+$, so the two
generators coincide. Put
\[
 y_0=\sqrt{R/2},\qquad p=e^{-R/2},
\]
and choose~$b>0$ with~$\mu_v([-b,b])=p$. Take independent initial
variables with
\begin{equation*}
 \sqrt\alpha X_0\sim N(y_0,1),\qquad V_0\sim\mu_v(\,\cdot\mid |v|\le b).
\end{equation*}
The translated spatial Gaussian has R\'enyi-2 divergence~$y_0^2=R/2$,
and the conditioned velocity law has R\'enyi-2 divergence~$-\log p=R/2$.
The R\'enyi-2 divergence is additive for product laws, so
$\Ren_2(\varrho_0\mid\mu)=R$. By symmetry, $\mathbb EV_0=0$, and the
truncation to~$[-b,b]$ gives~$\mathbb EV_0^2\le1$.

\emph{Speed bound.} Bounces and flips preserve~$|V|$. Hence the law
of~$|V_t|$ is the mixture of its initial law, with weight~$e^{-\gamma t}$,
and its equilibrium law, with weight~$1-e^{-\gamma t}$. Its density
relative to the equilibrium speed law is
\[
 1+e^{-\gamma t}\bigl(p^{-1}\mathbf1_{[0,b]}-1\bigr).
\]
Data processing for the map~$(x,v)\mapsto|v|$ gives
\begin{equation}\label{eq:sharp-speed-bound}
 \Ren_2(\varrho_t\mid\mu)\ge\log\left[1+e^{-2\gamma t}(e^{R/2}-1)\right].
\end{equation}
Since $\Ren_2(\varrho_t\mid\mu)\le1/2$, \eqref{eq:sharp-speed-bound}
and~$R\ge e^8$ imply
\begin{equation}\label{eq:sharp-speed-time}
 t\ge\frac1{2\gamma}\log\frac{e^{R/2}-1}{e^{1/2}-1}\ge\frac{R}{4\gamma}.
\end{equation}
The same mixture gives~$\mathbb EV_s^2=1+e^{-\gamma s}(\mathbb EV_0^2-1)\le1$
for all~$s\ge0$.

\emph{Return-time bound.} We measure the distance from the origin on a
logarithmic scale by~$\log(1+\alpha x^2)$, and add its transport
derivative divided by~$\gamma$, so that refreshment cancels that
derivative. Define
\begin{equation*}
 \Phi_\gamma(x,v)=\log(1+\alpha x^2)+\frac{2\alpha xv}{\gamma(1+\alpha x^2)}.
\end{equation*}
Write~$y=\sqrt\alpha x$. After this cancellation, the transport derivative
of the velocity term and the bounce term give
\begin{align}\label{eq:sharp-log-drift}
 \mc L\Phi_\gamma(x,v)
 &=\frac\alpha\gamma v^2\left[
 \frac{2(1-y^2)}{(1+y^2)^2}
 -\frac{4y^2}{1+y^2}\mathbf1_{\{yv>0\}}\right]\nonumber\\
 &\ge-\frac{4\alpha}{\gamma}v^2.
\end{align}
The inequality holds because adding~$4$ to the bracket and replacing the
indicator by~$1$ gives~$(6+2y^2)/(1+y^2)^2\ge0$. Dynkin's formula
and~$\mathbb EV_s^2\le1$ yield
\begin{equation}\label{eq:sharp-log-expectation}
 \mathbb E\Phi_\gamma(X_t,V_t)\ge\mathbb E\Phi_\gamma(X_0,V_0)-\frac{4\alpha t}{\gamma}.
\end{equation}
We first apply Dynkin's formula with localization on a finite
interval~$[0,T]$. Bounces preserve speed, and the expected sum of squared
refreshed velocities is~$\gamma T$, so
$\mathbb E\sup_{s\le T}|V_s|^2\le1+\gamma T$. The identity
$X_s=X_0+\int_0^sV_r\ud r$ then gives
$\mathbb E\sup_{s\le T}|X_s|^2<\infty$. Since~$\Phi_\gamma$ has at most
linear growth and~$|\mc L\Phi_\gamma|\le4\alpha v^2/\gamma$, these bounds
justify removing the localization.

It remains to compare the initial and terminal values
of~$\mathbb E\Phi_\gamma$. The velocity term of~$\Phi_\gamma$ has initial
mean zero. For~$Z\sim N(0,1)$, the
inequality~$1+y_0^2\le2[1+(y_0+Z)^2](1+Z^2)$ and Jensen's inequality give
\begin{equation*}
 \mathbb E\Phi_\gamma(X_0,V_0)\ge\log(1+R/2)-2\log2.
\end{equation*}
At equilibrium, Jensen's inequality gives~$\int\Phi_\gamma\ud\mu\le\log2$.
The Gaussian Poincar\'e inequality and $|2Z/(1+Z^2)|\le1$ bound the
variance of the logarithmic term by~$1$. The velocity term has
conditional mean zero given the position and variance at
most~$\alpha/\gamma^2$, so
\begin{equation*}
 \operatorname{Var}_\mu \Phi_\gamma\le1+\frac\alpha{\gamma^2}.
\end{equation*}
Put~$\zeta=\sqrt{e^{1/2}-1}<1$. Since~$\Ren_2(\varrho_t\mid\mu)\le1/2$, we
have~$\chi^2(\varrho_t\mid\mu)\le\zeta^2$, and Cauchy--Schwarz gives
\[
 \mathbb E\Phi_\gamma(X_t,V_t)\le\log2+\zeta\sqrt{1+\alpha/\gamma^2}.
\]
Combining this with~\eqref{eq:sharp-log-expectation} and
using~$\sqrt{1+\alpha/\gamma^2}\le1+\sqrt\alpha/\gamma$ yields
\begin{equation*}
 t+\frac{\zeta}{4\sqrt\alpha}\ge\frac\gamma{4\alpha}
 \left[\log(1+R/2)-3\log2-\zeta\right].
\end{equation*}
Since~$\log R\ge8$ and~$4\log2+1<4$, the bracket is at
least~$\tfrac12\log R$.

\emph{Conclusion.} Combining the last display
with~\eqref{eq:sharp-speed-time} gives
\begin{equation*}
 t+\frac1{4\sqrt\alpha}
 \ge\max\left\{\frac{R}{4\gamma},\frac{\gamma\log R}{8\alpha}\right\}
 \ge\frac1{\sqrt\alpha}\sqrt{\frac{R\log R}{32}}.
\end{equation*}
Hence $t\ge(16\sqrt\alpha)^{-1}\sqrt{R\log R}$ for~$R\ge e^8$. Since
$\gamma>0$ was arbitrary and the initial law does not depend on~$\gamma$,
this proves~\eqref{eq:initial-divergence-lower-bound}.
\end{proof}

\section{Lower bounds for the Zigzag sampler}
\label{sec:zz-obstruction}

The ZZ bound in Corollary~\ref{cor:mixing} depends on the initial
divergence through~$\CLSI^{-1/2}\sqrt{R\mathfrak{L}}$ and on the upper Hessian
bound through~$\CLSI^{-1}\sqrt\beta\,\mathfrak{L}$. This appendix shows that, in
dimension~$d\ge2$, neither dependence can be improved by more than a
logarithmic factor, in the worst case over Gaussian targets and initial
laws with bounded R\'enyi-2 divergence, uniformly in the refreshment rate.
The restriction $d\ge2$ is needed for the dependence on~$\beta$: in
dimension one, BPS and ZZ coincide, so the $\beta$-free entry
of~\eqref{eq:rate-bps} also bounds ZZ. The first example is
an isotropic Gaussian target with an initial law displaced along a
coordinate axis. The second is a family of two-dimensional Gaussian
targets confined in a direction oblique to the coordinate axes; it also
shows that strong convexity alone gives no uniform ZZ bound, in contrast
with BPS.

\begin{theorem}[Dependence on the initial divergence]
\label{thm:zz-finite-speed}
Let $d\ge1$ and $\alpha>0$, set $U(x)=\alpha|x|^2/2$ on~$\RR^d$, and
write $\mu=N(0,\alpha^{-1}I_d)\otimes\mu_v$, so that \eqref{eq:lsi} holds
with $\CLSI=\alpha$. For $R>0$, let
\[
 \varrho_0=N\bigl(\sqrt{R/\alpha}\,e_1,\alpha^{-1}I_d\bigr)\otimes\mu_v.
\]
Then $\Ren_2(\varrho_0\mid\mu)=R$, and the law~$\varrho_t$ of ZZ with
every refreshment rate~$\gamma>0$ satisfies, for all $t\ge0$,
\begin{equation}\label{eq:zz-finite-speed}
 \TV(\varrho_t,\mu)\ge\frac12\bigl(1-e^{-R/8}\bigr)
 -2\sqrt{\frac{2\alpha}{\pi R}}\;t.
\end{equation}
In particular, if $R\ge16$, then $\TV(\varrho_t,\mu)>1/4$ for
$0\le t\le\frac1{10}\sqrt{R/\alpha}$.
\end{theorem}

\begin{proof}
The position marginals of~$\varrho_0$ and~$\mu$ are Gaussians with the
same covariance~$\alpha^{-1}I_d$ and with means $\sqrt{R/\alpha}\,e_1$
and~$0$, and the velocity marginals coincide. For two Gaussians with a
common covariance~$\Sigma$ and means $m$ and~$0$, the R\'enyi-2
divergence equals $m^{\mathsf T}\Sigma^{-1}m$, so
$\Ren_2(\varrho_0\mid\mu)=\alpha\cdot R/\alpha=R$.

Let $(X_t,V_t)$ be the ZZ process with initial law~$\varrho_0$. A flip
changes the sign of one velocity coordinate and preserves its absolute
value, and a refreshment redraws~$V$ from~$\mu_v$, independently of the
past. Hence $|V_{s,1}|$ equals $|V_{0,1}|$ if no refreshment occurs
in~$[0,s]$, and otherwise the absolute first coordinate of the last
refreshed velocity. Both have the law of~$|Z|$ with $Z\sim N(0,1)$ and
are independent of the refreshment times, so
$\mathbb E|V_{s,1}|=\sqrt{2/\pi}$ for every $s\ge0$, whatever the
positions. Write $Y_t=\sqrt\alpha\,X_{t,1}$ and
$D_t=\sqrt\alpha\int_0^t|V_{s,1}|\ud s$, so that $Y_t\ge Y_0-D_t$ and
$\mathbb ED_t=\sqrt{2\alpha/\pi}\;t$.

Consider the event $A=\{x\in\RR^d:\sqrt\alpha\,x_1\ge\sqrt R/2\}$. Under
$\mu$, $\sqrt\alpha\,x_1\sim N(0,1)$, and the Gaussian tail bound gives
$\mu(A\times\RR^d)\le\frac12e^{-R/8}$. Under $\varrho_t$, since
$Y_0\sim N(\sqrt R,1)$ has median~$\sqrt R$, Markov's inequality gives
\[
 \varrho_t(A\times\RR^d)
 \ge\mathbb P\bigl(Y_0\ge\sqrt R,\;D_t\le\sqrt R/2\bigr)
 \ge\frac12-\frac{2\,\mathbb ED_t}{\sqrt R}
 =\frac12-2\sqrt{\frac{2\alpha}{\pi R}}\;t.
\]
The difference of the two probabilities of~$A\times\RR^d$ is a lower
bound for the total variation distance, which
proves~\eqref{eq:zz-finite-speed}. For $R\ge16$ and
$t\le\frac1{10}\sqrt{R/\alpha}$, the right-hand side is at least
$\tfrac12(1-e^{-2})-\tfrac15\sqrt{2/\pi}$, which exceeds~$1/4$.
\end{proof}

For this target, Corollary~\ref{cor:mixing} with $\beta=\CLSI=\alpha$
gives $T^{\mathrm{TV}}_{1/4}\le C\alpha^{-1/2}\sqrt{R\log R}$ for a
universal constant~$C$, so the lower bound
$\frac1{10}\sqrt{R/\alpha}$ is sharp up to the factor~$\sqrt{\log R}$.
Theorem~\ref{thm:initial-divergence-sharpness} shows that in R\'enyi-2
divergence this factor cannot be removed. In dimension one, the BPS and
ZZ generators coincide, so Theorem~\ref{thm:zz-finite-speed} also applies
to BPS when $d=1$. In higher dimensions, bounces preserve only the
speed~$|V|$, whose mean under~$\mu_v$ is at most~$\sqrt d$, and the same
argument gives for BPS only a lower bound of order~$\sqrt{R/(\alpha d)}$.

\begin{theorem}[Dependence on the upper Hessian bound]
\label{thm:zz-obstruction}
For $\kappa\ge1$, introduce the orthogonal coordinates
\[
 s=\frac{x_1+x_2}{\sqrt2},\qquad r=\frac{x_1-x_2}{\sqrt2},
\]
and let
\begin{equation}\label{eq:zz-obstruction-potential}
 U_\kappa(x)=\frac{s^2}{2}+\frac{\kappa r^2}{2},\qquad
 \mu_\kappa(\ud x,\ud v)\propto e^{-U_\kappa(x)}\ud x\,\mu_v(\ud v),
 \quad \mu_v=N(0,I_2).
\end{equation}
These potentials satisfy $I_2\preceq\nabla^2U_\kappa\preceq \kappa I_2$, and the
position marginals of~$\mu_\kappa$ satisfy \eqref{eq:lsi} with $\CLSI=1$.
Let $\varrho_t^{\kappa,\gamma}$ be the law of ZZ with
potential $U_\kappa$, refreshment rate $\gamma>0$, and initial law satisfying
\[
 f_0=\frac{\ud\varrho_0^{\kappa,\gamma}}{\ud\mu_\kappa}-1=\frac12\sin s,
\]
and write $f_t=\ud\varrho_t^{\kappa,\gamma}/\ud\mu_\kappa-1$.
The initial R\'enyi-2 divergence and total variation distance are
\[
 R_*=\log\left(1+\frac{1-e^{-2}}8\right)>0,
 \qquad \theta_*=\frac14\mathbb E|\sin Z|>\frac1{10},
\]
where $Z\sim N(0,1)$. For every $\gamma>0$ and $t\ge0$,
\begin{equation}\label{eq:zz-obstruction-quantitative}
 \norm{f_t-f_0}_{L^2(\mu_\kappa)}^2
 \le\frac{t}{\sqrt\kappa}\bigl[2+\log(1+t\sqrt\kappa)\bigr].
\end{equation}
Consequently, for every $\gamma>0$,
\begin{equation}\label{eq:zz-obstruction-times}
 \TV(\varrho_t^{\kappa,\gamma},\mu_\kappa)\ge\frac{\theta_*}2
 \quad\text{and}\quad
 \Ren_2(\varrho_t^{\kappa,\gamma}\mid\mu_\kappa)\ge\frac{R_*}2
 \qquad\text{for }0\le t\le\frac{\sqrt\kappa}{500\log(e\kappa)}.
\end{equation}
\end{theorem}

The time scale in~\eqref{eq:zz-obstruction-times} tends to infinity as~$\kappa\to\infty$,
so no bound on the ZZ mixing time that is uniform over $\kappa\ge1$ can
hold for this family, although $d$, $\CLSI$, the strong-convexity
constant, and the initial divergence are fixed.

\begin{proof}
The Hessian eigenvalues in the $(s,r)$ coordinates are $1$ and~$\kappa$,
so the curvature assertions and the LSI with $\CLSI=1$ follow from strong
convexity and the Bakry--\'Emery criterion~\cite{BakryEmery1985}.
Under~$\mu_\kappa$, $s$ and $\sqrt\kappa\,r$ are independent standard
Gaussians, which gives $\norm{f_0}^2=(1-e^{-2})/8$ and the values
of~$R_*$ and~$\theta_*$. The bound $\theta_*>1/10$ follows from
$\mathbb E|\sin Z|\ge\mathbb E\sin^2Z=(1-e^{-2})/2$. We
prove~\eqref{eq:zz-obstruction-quantitative} by showing that the initial
spatial density changes little under the adjoint semigroup. The change is
controlled by a resolvent of the symmetric part of the generator, whose
spectral parameter is matched to the time horizon.

Fix $\kappa\ge1$ and $\gamma>0$, and write $\mc P_t=e^{t\mc L}$ for the
ZZ semigroup. All norms and
inner products below are in $L^2(\mu_\kappa)$.
On test functions, the symmetric part $\mc S=(\mc L+\mc L^*)/2$ of the
generator acts only on velocity:
\begin{equation}\label{eq:zz-obstruction-symmetric}
 \mc Sh=\frac12\sum_{i=1}^2|v_i\partial_{x_i}U_\kappa|
       [h(x,F_iv)-h(x,v)]+\gamma(\Pi_v h-h).
\end{equation}
Its nonnegative quadratic form is
\begin{equation}\label{eq:zz-obstruction-form}
 \mc E(h)=\frac14\sum_{i=1}^2
 \int |v_i\partial_{x_i}U_\kappa|[h(x,F_iv)-h(x,v)]^2\ud\mu_\kappa
 +\gamma\norm{(I-\Pi_v)h}^2.
\end{equation}
Write $\mc E(h,k)$ for the associated bilinear form.
The ZZ dissipation identity in Lemma~\ref{lem:velocity-dissipation},
which holds without the zero-mean assumption, gives
$\mc E(h)=-\ip{h}{\mc Lh}$ for test functions~$h$: the change of
variables $v\mapsto F_iv$ preserves $\mu_\kappa$ and the squared jump and
reverses the sign of $v_i\partial_{x_i}U_\kappa$, so the one-sided rates
in the lemma may be replaced by the symmetric weights
in~\eqref{eq:zz-obstruction-form}. Graph-norm approximation extends this
identity to $h\in D(\mc L)$. Consequently,
\begin{equation}\label{eq:zz-obstruction-semigroup-energy}
 \int_0^t\mc E(\mc P_\tau h)\ud\tau
 =\frac12\bigl(\norm h^2-\norm{\mc P_t h}^2\bigr)
 \le\frac12\norm h^2,
 \qquad h\in D(\mc L).
\end{equation}

\emph{A parametrized resolvent.} The velocity observable
$(v_1+v_2)/\sqrt2$ has an explicit resolvent for~$\mc S$. For $a>0$,
define
\begin{equation}\label{eq:zz-obstruction-test}
 \xi_a(x,v)=\frac1{\sqrt2}\sum_{i=1}^2
 \frac{v_i}{a+\gamma+|v_i\partial_{x_i}U_\kappa(x)|}.
\end{equation}
Each summand is odd in its own velocity coordinate and independent of
the other. Thus $\Pi_v\xi_a=0$, the flip~$F_i$ reverses the sign of the
$i$th summand and fixes the other, and
\begin{equation}\label{eq:zz-obstruction-resolvent}
 (a-\mc S)\xi_a=\frac{v_1+v_2}{\sqrt2}.
\end{equation}
Set
\begin{equation}\label{eq:zz-obstruction-b}
 b_\kappa(a)=\frac12\sum_{i=1}^2
 \int\frac{v_i^2}{a+|v_i\partial_{x_i}U_\kappa(x)|}\ud\mu_\kappa.
\end{equation}
Parity eliminates mixed velocity terms in the norm and in the form, and
direct substitution gives
\begin{equation}\label{eq:zz-obstruction-energy}
 a\norm{\xi_a}^2+\mc E(\xi_a)
 =\frac12\sum_{i=1}^2
 \int\frac{v_i^2}{a+\gamma+|v_i\partial_{x_i}U_\kappa(x)|}\ud\mu_\kappa
 \le b_\kappa(a).
\end{equation}

\emph{Almost invariance of~$f_0$.} As a function of~$x$ alone with
bounded derivative, $f_0$ belongs to $D(\mc L)\cap D(\mc L^*)$ and
satisfies $\mc Lf_0=v\cdot\nabla_xf_0=-\mc L^*f_0$.
Put $\omega_a=\tfrac12\cos(s)\xi_a$.
Since $\mc S$ acts only on velocity and $\cos^2s\le1$,
\[
 (a-\mc S)\omega_a
 =\frac{\cos s}{2}\frac{v_1+v_2}{\sqrt2}
 =\mc Lf_0=-\mc L^*f_0,
 \qquad
 a\norm{\omega_a}^2+\mc E(\omega_a)\le\frac{b_\kappa(a)}{4}.
\]
For test functions $h$, changing variables under each flip gives
\[
 \ip{(a-\mc S)\omega_a}{h}=a\ip{\omega_a}{h}+\mc E(\omega_a,h).
\]
Graph-norm approximation extends this identity to $h\in D(\mc L)$; since
it involves only $\omega_a$ and its velocity differences, $\omega_a$ need
not belong to $D(\mc L)$.

Duhamel's formula and
\eqref{eq:zz-obstruction-semigroup-energy} imply, for $h\in D(\mc L)$,
\begin{align*}
 |\ip{f_t-f_0}{h}|
 &=\left|\int_0^t
   \left[a\ip{\omega_a}{\mc P_\tau h}
             +\mc E(\omega_a,\mc P_\tau h)\right]\ud\tau\right|\\
 &\le\norm h\left[at\norm{\omega_a}+
                    \sqrt{\frac t2\mc E(\omega_a)}\right].
\end{align*}
Here the first term uses $L^2$ contractivity, and the second uses
Cauchy--Schwarz in time and in the form.
Taking the supremum over the dense domain $D(\mc L)$ and then applying
Cauchy--Schwarz in $\RR^2$ gives
\begin{equation}\label{eq:zz-obstruction-density}
 \norm{f_t-f_0}^2
 \le\left(at^2+\frac t2\right)
          [a\norm{\omega_a}^2+\mc E(\omega_a)]
 \le\frac14\left(at^2+\frac t2\right)b_\kappa(a).
\end{equation}
The right-hand side is independent of $\gamma$.

\emph{Bound on~$b_\kappa(a)$.} Put $y=\sqrt\kappa\,r$, so that $s$
and $y$ are independent standard Gaussians under $\mu_\kappa$ for every
$\kappa$. Moreover,
\[
 \partial_{x_1}U_\kappa
 =\frac{s+\sqrt\kappa\,y}{\sqrt2},\qquad
 \partial_{x_2}U_\kappa
 =\frac{s-\sqrt\kappa\,y}{\sqrt2}.
\]
Each partial derivative is a centered Gaussian with
variance~$\sigma_\kappa^2=(1+\kappa)/2$, independent of~$v$. Hence
\[
 b_\kappa(a)=\mathbb E\,\frac{V^2}{a+\sigma_\kappa|VZ|},
\]
where $V$ and $Z$ are independent standard Gaussians. For $c>0$, we
split the integral over~$Z$ at $|Z|=1$. On $\{|Z|\le1\}$ the Gaussian
density is at most $(2\pi)^{-1/2}<1/2$, and on $\{|Z|>1\}$ we
bound $a+c|Z|$ below by~$c$. Since $\mathbb P(|Z|>1)<1$,
\[
 \mathbb E\,\frac1{a+c|Z|}
 \le\frac1c\log\Bigl(1+\frac ca\Bigr)+\frac1c
 =\frac1c\Bigl[1+\log\Bigl(1+\frac ca\Bigr)\Bigr].
\]
We condition on~$V$ and apply this bound with $c=\sigma_\kappa|V|$.
Then $\log(1+xy)\le\log(1+x)+\log(1+y)$ for $x,y\ge0$,
$\mathbb E|V|\le1$, and $|V|\log(1+|V|)\le V^2$ give
\begin{equation}\label{eq:zz-obstruction-b-bound}
 b_\kappa(a)
 \le\frac1{\sigma_\kappa}\,
 \mathbb E\Bigl[|V|\Bigl(1+\log\Bigl(1+\frac{\sigma_\kappa|V|}a\Bigr)\Bigr)\Bigr]
 \le\frac{2+\log(1+\sigma_\kappa/a)}{\sigma_\kappa}.
\end{equation}

\emph{Conclusion.} The case $t=0$ is trivial. For $t>0$, take $a=1/t$
in \eqref{eq:zz-obstruction-density}
and~\eqref{eq:zz-obstruction-b-bound}. Since
$\sqrt{\kappa/2}\le\sigma_\kappa\le\sqrt\kappa$ for $\kappa\ge1$ and
$3\sqrt2/8<1$,
\[
 \norm{f_t-f_0}^2
 \le\frac{3t}{8\sigma_\kappa}\bigl[2+\log(1+\sigma_\kappa t)\bigr]
 \le\frac{t}{\sqrt\kappa}\bigl[2+\log(1+t\sqrt\kappa)\bigr],
\]
which is \eqref{eq:zz-obstruction-quantitative}.

For \eqref{eq:zz-obstruction-times}, let
$0\le t\le\sqrt\kappa/[500\log(e\kappa)]$. Then $t\sqrt\kappa\le\kappa$,
so $\log(1+t\sqrt\kappa)\le\log(e\kappa)$,
and~\eqref{eq:zz-obstruction-quantitative} gives
\[
 \norm{f_t-f_0}^2\le\frac{2+\log(e\kappa)}{500\log(e\kappa)}\le\frac3{500},
 \qquad\text{so}\qquad \norm{f_t-f_0}<0.078.
\]
Since $\norm{f_t-f_0}_{L^1(\mu_\kappa)}\le\norm{f_t-f_0}$
and~$\theta_*>1/10$,
\[
 \TV(\varrho_t^{\kappa,\gamma},\mu_\kappa)
 \ge\theta_*-\frac12\norm{f_t-f_0}
 >\theta_*-0.039>\frac{\theta_*}2.
\]
For the R\'enyi-2 divergence, $\norm{f_0}^2=(1-e^{-2})/8>0.108$ gives
$\norm{f_t}\ge\norm{f_0}-\norm{f_t-f_0}>0.328-0.078=0.25$. Hence
$\Ren_2(\varrho_t^{\kappa,\gamma}\mid\mu_\kappa)=\log(1+\norm{f_t}^2)
>\log(1.0625)>0.06$, while $R_*<\log(1.109)<0.104$. This
proves~\eqref{eq:zz-obstruction-times}.
\end{proof}

The two examples combine into a lower bound with the structure of the
ZZ entry in~\eqref{eq:optimized-mixing-times}.

\begin{corollary}[Sharpness of the Zigzag bound]\label{cor:zz-sharpness}
Let $d\ge2$, $0<\CLSI\le\beta<\infty$, $R\ge16$, and $0<\delta\le1/20$.
There exist a potential~$U$ on~$\RR^d$ with
$\CLSI I_d\preceq\nabla^2U\preceq\beta I_d$, so that \eqref{eq:lsi}
holds with constant~$\CLSI$, and an initial law~$\varrho_0$ with
$\Ren_2(\varrho_0\mid\mu)\le R$, such that ZZ with every refreshment
rate~$\gamma>0$ satisfies
\begin{equation}\label{eq:zz-sharpness}
 T_\delta^{\mathrm{TV}}
 \ge\frac1{20}\sqrt{\frac R{\CLSI}}
 +\frac1{1000}\,\frac{\sqrt\beta}{\CLSI\log(e\beta/\CLSI)}.
\end{equation}
\end{corollary}

\begin{proof}
We construct two examples, the first with
$T_\delta^{\mathrm{TV}}\ge t_1=\frac1{10}\sqrt{R/\CLSI}$ and the second with
$T_\delta^{\mathrm{TV}}\ge t_2=\frac{\sqrt\beta}{500\,\CLSI\log(e\beta/\CLSI)}$,
each for every refreshment rate. We take the first example if $t_1\ge t_2$
and the second otherwise. This choice depends only on $d$, $\CLSI$,
$\beta$, and~$R$, not on~$\gamma$. The chosen example satisfies
$T_\delta^{\mathrm{TV}}\ge\max\{t_1,t_2\}\ge(t_1+t_2)/2$ for every~$\gamma$,
which is~\eqref{eq:zz-sharpness}.

For the first, take $U(x)=\CLSI|x|^2/2$ and the initial law of
Theorem~\ref{thm:zz-finite-speed} with $\alpha=\CLSI$. The Hessian
is~$\CLSI I_d\preceq\beta I_d$, and $\Ren_2(\varrho_0\mid\mu)=R$. Since
$R\ge16$ and $\delta\le1/4$, Theorem~\ref{thm:zz-finite-speed} gives
$\TV(\varrho_t,\mu)>1/4\ge\delta$ for $t\le\frac1{10}\sqrt{R/\CLSI}$.

For the second, put $\kappa=\beta/\CLSI\ge1$. We first rescale. For
$c>0$ and a potential~$W$ on~$\RR^2$, let $W_c(x)=W(\sqrt c\,x)$. If
$(X_t,V_t)$ is ZZ for~$W$ with refreshment rate~$\gamma$, then
$(X_{\sqrt c\,t}/\sqrt c,\,V_{\sqrt c\,t})$ is ZZ for~$W_c$ with
refreshment rate~$\sqrt c\,\gamma$: the position has derivative
$V_{\sqrt c\,t}$, the flip rate of coordinate~$i$ at time~$t$ is
$\sqrt c\,(V_{\sqrt c\,t,i}\,\partial_iW(X_{\sqrt c\,t}))_+
=(V_{\sqrt c\,t,i}\,\partial_iW_c(X_{\sqrt c\,t}/\sqrt c))_+$, and the
refreshment rate is multiplied by~$\sqrt c$. The map
$(x,v)\mapsto(x/\sqrt c,v)$ sends the equilibrium of~$W$ to that
of~$W_c$ and preserves total variation distances and R\'enyi
divergences. Apply this with $W=U_\kappa$ and $c=\CLSI$. Since
$\gamma\mapsto\sqrt{\CLSI}\,\gamma$ is a bijection of~$(0,\infty)$,
Theorem~\ref{thm:zz-obstruction} shows that ZZ for~$W_{\CLSI}$, with any
refreshment rate and with the rescaled initial law, has total variation
distance at least $\theta_*/2>1/20\ge\delta$ from equilibrium for
$0\le t\le\sqrt\kappa/[500\sqrt{\CLSI}\log(e\kappa)]$, and
$\sqrt\kappa/\sqrt{\CLSI}=\sqrt\beta/\CLSI$.

Now define, for $x\in\RR^d$,
\[
 U(x)=U_\kappa\bigl(\sqrt{\CLSI}\,(x_1,x_2)\bigr)+\frac{\CLSI}2\sum_{i=3}^dx_i^2,
\]
whose Hessian has eigenvalues $\CLSI$ and~$\beta$, and let $\varrho_0$
be the product of the rescaled initial law in the variables
$(x_1,x_2,v_1,v_2)$ with the equilibrium law of the remaining
coordinates. R\'enyi-2 divergence is additive for products, so
$\Ren_2(\varrho_0\mid\mu)=R_*<1\le R$. For $i\in\{1,2\}$, the partial
derivative $\partial_{x_i}U$ depends only on~$(x_1,x_2)$, and a
refreshment redraws $(v_1,v_2)$ from~$N(0,I_2)$. Hence
$(X_{t,1},X_{t,2},V_{t,1},V_{t,2})$ is a ZZ process for~$W_{\CLSI}$ with
the same refreshment rate, and its law is the marginal of~$\varrho_t$.
Marginalization does not increase the total variation distance, so
$\TV(\varrho_t,\mu)>\delta$ for
$0\le t\le\sqrt\beta/[500\,\CLSI\log(e\beta/\CLSI)]$.
\end{proof}

For the parameters of Corollary~\ref{cor:zz-sharpness},
Corollary~\ref{cor:mixing}, with the refreshment rate chosen as in that corollary, gives
\[
 T_\delta^{\mathrm{TV}}
 \le\frac{6000\,\mathfrak{L}}{\sqrt{\CLSI}}\sqrt{1+\frac R{\mathfrak{L}}+\frac\beta{\CLSI}}
 \le\frac{6000}{\sqrt{\CLSI}}
 \Bigl(\sqrt{R\mathfrak{L}}+2\mathfrak{L}\sqrt{\beta/\CLSI}\Bigr),
 \qquad \mathfrak{L}=\log\frac{3R}{2\delta^2},
\]
where we used $\mathfrak{L}\le\mathfrak{L}\sqrt{\beta/\CLSI}$. The ratio of this upper
bound to the lower bound~\eqref{eq:zz-sharpness} is at most a universal
constant times~$\mathfrak{L}\log(e\beta/\CLSI)$. Hence, for $d\ge2$, no bound on
the ZZ mixing time in terms of $d$, $\CLSI$, $\beta$, $R$, and~$\delta$
alone can improve on Corollary~\ref{cor:mixing} by more than a factor
logarithmic in~$R$, $1/\delta$, and~$\beta/\CLSI$. Letting
$\beta\to\infty$ with the other parameters fixed shows that strong
convexity alone provides no uniform mixing bound for ZZ.

The lower bounds also hold for the position marginal, for which
Corollary~\ref{cor:mixing} gives the same upper bound by data processing.
In Theorem~\ref{thm:zz-finite-speed}, the event~$A\times\RR^d$ depends
only on the position. In Theorem~\ref{thm:zz-obstruction}, the initial
density~$f_0$ depends only on~$x$, so
$\norm{\Pi_vf_t-f_0}=\norm{\Pi_v(f_t-f_0)}\le\norm{f_t-f_0}$, and the
argument for~\eqref{eq:zz-obstruction-times} applies to~$\Pi_vf_t$. The
obstruction is therefore not an effect of slow velocity equilibration
alone.

For the family of Theorem~\ref{thm:zz-obstruction}, the $\beta$-free
entry of~\eqref{eq:rate-bps} gives a BPS bound uniform in~$\kappa$, since
$\CLSI=1$ and $d=2$ are fixed. BPS reflects the
velocity across~$\nabla U_\kappa$ rather than reversing coordinates. The
valley component of the reflected velocity changes by an amount
proportional to~$s/|\nabla U_\kappa|$, which is small wherever
$\kappa|r|\gg|s|$, that is, for most of~$\mu_\kappa$ when $\kappa$ is
large. ZZ flips instead randomize the valley component at a rate that
grows with~$\kappa$. Theorem~\ref{thm:zz-obstruction} thus identifies a
genuine difference between the two dynamics. The lower bounds of this
appendix are worst-case statements over the stated classes of targets
and initial laws. They do not assert failure of convergence for a fixed
target, and they do not concern the Gaussian cold start
specifically.

\bibliography{references}

\end{document}